\input ./style/arxiv-general.cfg
\documentclass[MSNbibl,number,citesort,seceqn,dvips]{arxbj}
\makeatletter
   \@ifpackageloaded{graphicx}{}{\usepackage{graphicx}}
\makeatother
\usepackage{mathbh}
\usepackage{stmaryrd}

\volume{22}
\issue{1}
\pubyear{2016}
\firstpage{345}
\lastpage{375}
\doi{10.3150/14-BEJ661} 
\docsubty{FLA}

\makeatletter
\newtheorem{theorem}{Theorem}[section]
\newtheorem{cor}[theorem]{Corollary}
\newcommand{\ds}{\displaystyle}
\newcommand{\rrvert}{\vert}
\newcommand{\llvert}{\vert}
\newremark{exa}[theorem]{Example}
\newremark{rem}[theorem]{Remark}
\newtheorem{prop}[theorem]{Proposition}
\makeatother

\begin{document}
\begin{frontmatter}

\title{Central limit theorems for long range dependent spatial linear processes}
\runtitle{CLT under spatial long memory}

\begin{aug}
\author[A]{\inits{S.N.}\fnms{S.N.}~\snm{Lahiri}\corref{}\thanksref{A}\ead[label=e1]{snlahiri@stat.ncsu.edu}}
\and
\author[B]{\inits{P.M.}\fnms{Peter M.}~\snm{Robinson}\thanksref{B}\ead[label=e2]{p.m.robinson@lse.ac.uk}}
\address[A]{Department of Statistics,
North Carolina State University,
2311 Stinson Dr,
Raleigh, NC 27695-8203,
USA. \printead{e1}}
\address[B]{Department of Economics,
London School of Economics,
Houghton Street,
London WC2A 2AE,
UK.\\ \printead{e2}}
\end{aug}

\received{\smonth{4} \syear{2013}}
\revised{\smonth{5} \syear{2014}}

%
\begin{abstract}
Central limit theorems are established for the sum, over a spatial region,
of observations from a linear process on a $d$-dimensional lattice. This
region need not be rectangular, but can be irregularly-shaped. Separate
results are established for the cases of positive strong dependence,
short range
dependence, and negative dependence. We provide approximations to
asymptotic variances that reveal differential rates of convergence under
the three types of dependence. Further, in contrast to the
one dimensional (i.e., the time series) case,
it is shown that the
form of the asymptotic variance in dimensions $d>1$
critically depends on the geometry of the sampling region
under positive strong dependence and under negative dependence
and that there can be non-trivial edge-effects
under negative dependence for $d>1$. Precise conditions
for the presence of edge effects are also given.
\end{abstract}

%
\end{frontmatter}


\section{Introduction}\label{sec1}

The presence of long range dependence in spatial data has been noted in
various empirical studies but a suitable formulation and systematic
study of such spatial processes
is lacking. For example, the ``law of environmental variation''
of Fairfield~Smith \cite{Fai38}, based on ``Agricultural Field
Trials'' data,
posits that the covariance function of the
yield in the plane decays as the inverse of the Euclidean
distance. Thus, the covariance functions of such spatial
processes are not absolutely summable and may exhibit long
range dependence. More recently, the effect of spatial long range
dependence has been noted in Atmospheric sciences
(cf. Kashyap and Lapsa \cite{KasLap84}, Gneiting \cite{Gne00}), Economics (Leonenko and Taufer \cite{LeoTau13}),
Oceanography (cf. Percival \textit{et al.} \cite{Peretal08}) and
Solid State Physics (cf. Carlos-Davila \textit{et al.} \cite{CarMejMor85}),
among others. See also Lavancier \cite{Lav06} for some
specific examples and other applications of
spatial long range dependence.
The traditional approach of quantifying
spatial dependence through various notions of mixing
is inadequate for dealing with long range dependence.
In this paper, we consider a class of stationary spatial linear
processes that allow for long range dependence, as well
as the properties of short range- and negative-dependence,
and establish central limit theorems for the sum
over the entire
range of such dependence.

The dependence structure of a real-valued stationary process on a
$d$-dimensional spatial lattice can be non-parametrically modeled by
the linear process
%
\begin{equation}\label{1.1}
Z(\mathbf{i})=\mu+\sum_{\mathbf{j}\in\mathbb{Z}^{d}} \alpha ( \mathbf{i}-\mathbf{j})
\varepsilon(\mathbf{j}), \qquad \mathbf{i}\in \mathbb{Z}^{d},
\end{equation}
where the collection of real numbers  $ \{ \alpha ( \mathbf
{i}),\mathbf{i}\in\mathbb{Z}^{d} \} $ satisfies
%
\begin{equation}\label{1.2}
\sum_{\mathbf{i}\in\mathbb{Z}^{d}} \alpha ( \mathbf{i} )
^{2}<\infty,
\end{equation}
and $ \{ \varepsilon(\mathbf{i}), \mathbf{i}\in\mathbb{Z}^{d}\} $
is a collection of independent homoscedastic random variables with zero mean
and finite variance. If the $\varepsilon(\mathbf{i})$ are only
uncorrelated, (\ref{1.1}) and (\ref{1.2}) represent the class of purely non-deterministic
processes on $\mathbb{Z}^{d}$, for which $Z(\mathbf{i})$ has a more
parsimonious ``half-plane'' representation (see, e.g., Whittle \cite{Whi54}), to
generalize the one-sided Wold representation in the time series case $d=1$.
 However, we impose independence in order to establish central limit
theorems (CLTs) for
\[
S_{n}=\sum_{\mathbf{i}\in\mathcal{D}_{n}} Z(\mathbf{i})
\]
as $n\rightarrow\infty$.  Here it is supposed that we observe
$Z(\mathbf{i})$ within a spatial region $R_{n}\subset\mathbb{R}^{d}$ (to be described
in detail subsequently) whose volume is regarded as increasing with the
integer $n\geq 1$, the data sites being given by
\[
\mathcal{D}_{n}=R_{n}\cap\mathbb{Z}^{d}.
\]

Extending time series notions (cf. Robinson \cite{Rob97}), we consider three
different sub-classes of~(\ref{1.1}), (\ref{1.2}), broadly described as
\begin{eqnarray*}
 \mbox{negatively dependent  (ND): }\qquad\sum_{\mathbf{i}\in\mathbb{Z}^{d}}
\bigl\llvert \alpha ( \mathbf{i} ) \bigr\rrvert &<& \infty, \qquad\sum
_{\mathbf{i}\in\mathbb{Z}^{d}} \alpha ( \mathbf {i} ) =0,
\\
\mbox{short-range dependent (SRD): }\qquad \sum_{\mathbf{i}\in\mathbb{Z}^{d}}
\bigl\llvert \alpha ( \mathbf{i} ) \bigr\rrvert &<& \infty, \qquad\sum
_{\mathbf{i}\in\mathbb{Z}^{d}} \alpha ( \mathbf {i} ) \neq0,
\\
\mbox{positively strongly dependent (PSD): }\qquad\sum
_{\mathbf{i}\in\mathbb{Z}^{d}} \bigl\llvert \alpha ( \mathbf{i} ) \bigr\rrvert  &=& \infty,
\end{eqnarray*}
though our results rest on conditions that respectively, imply these. Denoting by $f (\mathcal{\lambda})$ the spectral
density of $Z(\mathbf{i})$, for ND processes $f(0) =0$, for SRD
processes $f(0)\in(0,\infty )$, and for PSD
processes $f$
may diverge at frequency $0$.  The three sub-classes are also
associated with
different rates of increase of
\[
\sigma_{n}^{2}=\operatorname{Var}(S_{n}),
\]
where we assume
%
\begin{equation}\label{1.3}
\sigma_{n}^{2}\rightarrow\infty\qquad\mbox{as }n\rightarrow
\infty,
\end{equation}
which is a necessary condition for a CLT.  Define $N_{n}=\llvert
\mathcal{D}_{n}\rrvert$,  where $\llvert  B\rrvert$
denotes the
size (i.e., number of elements) of a finite set $B$, so $N_{n}$ denotes
sample size.  Under additional conditions, we have
\begin{eqnarray*}
N_{n}^{-1}\sigma_{n}^{2} &
\rightarrow & 0 \qquad \mbox{as }n\rightarrow \infty\mbox{ when }Z(\mathbf{i})\mbox{ is ND},
\\
N_{n}^{-1}\sigma_{n}^{2} & \rightarrow & \sigma_{0}^{2}\in ( 0,\infty )\qquad  \mbox{as $n
\rightarrow\infty$ when $Z(\mathbf{i})$ is SRD},
\\
N_{n}^{-1}\sigma_{n}^{2} &
\rightarrow & \infty\qquad\mbox{as $n\rightarrow\infty$ when $Z(\mathbf{i})$ is PSD}.
\end{eqnarray*}
The condition of PSD is also referred to as
\textit{long range dependence} in the literature.

Under (\ref{1.1}), (\ref{1.2}) and (\ref{1.3}), and for $d=1$ with $\mathcal{D}_{n}=(1,2,\ldots,n)$,
Ibragimov and Linnik \cite{IbrLin71}, pages 359--360, established that
\begin{equation}\label{1.4}
[ S_{n}-ES_{n} ] /\sigma_{n}
\mathop{\rightarrow}^{d}\mathcal{N}(0,1)\qquad\mbox{as }n
\rightarrow\infty.
\end{equation}
Their main achievement was to allow arbitrarily slowly increasing
$\sigma_{n}^{2}$, in particular to cover all ND $Z(\mathbf{i})$, as well as
SRD and
PSD ones.
Again for $d=1$,  Hannan \cite{Han79} relaxed independence of the
$\varepsilon(\mathbf{i})$ to a martingale difference assumption, but only
covered SRD and PSD $Z(\mathbf{i})$.  On the other hand, Rosenblatt \cite{Ros61},
Taqqu \cite{Taq74} and others established non-central limit theorems when $d=1$
and $Z(\mathbf{i})$ does not satisfy (\ref{1.1}) but is a non-linear
function of a
PSD Gaussian process, and more generally.  Mention must also be made
of the
many CLTs for $d=1$, where (\ref{1.1}) is replaced by mixing conditions, following
Rosenblatt \cite{Ros56}, implying $Z(\mathbf{i})$ is SRD, and extended to $d>1$
by a number of authors; see, for example, Bolthausen \cite{Bol82}, Doukhan \cite{Dou94},
Guyon \cite{Guy95}, the latter two authors
also discussing the mixing properties of linear processes. CLTs
for SRD spatial processes over irregular sampling regions
are given by Lahiri \cite{Lah03} and El~Machkouri, Voln{\'y} and Wu \cite{ElMVolWu13},
allowing more general processes than linear fields.
But relatively less attention
has been paid to PSD processes with $d>1$, under either linear or other
assumptions. For rectangular regions, CLTs and invariance principles
for PSD spatial linear
processes have been proved by Lavancier \cite{Lav07} and for fractional Brownian
sheets 
by Wang \cite{Wan13}. Dobrushin and Major \cite{DobMaj79} and Surgailis \cite{Sur82}
proved central- and non-central limit theorems for functionals of PSD
Gaussian processes
and
for functionals of
PSD linear fields,
respectively.  There is also a small body of
literature on statistical inference on the mean and covariance parameters
of PSD spatial processes;
see the papers by Boissy \textit{et al.} \cite{Boietal05}, Beran \textit{et al.} \cite{BerGhoSch09}
and Wang and Cai \cite{WanCai10}
and
the monographs by Ivanov and Leonenko \cite{IvaLeo89} and Bertail \textit{et al.} \cite{bert06},
and the references therein. We know of no spatial work
under ND.

A major limitation of the existing work on spatial PSD processes is
that it
deals exclusively with rectangular spatial sampling regions. In
contrast to
the temporal case, in most practical applications spatial sampling regions
are non-rectangular, and possibly of a non-standard shape (cf. Cressie \cite{Cre93},
Lahiri \textit{et al.} \cite{Lahetal99}). As a result,
existing results are of limited use. The present paper attempts to fill the
gap by introducing a general framework for studying linear spatial processes
over sampling regions of non-standard shapes. The main results of the paper
establish separate CLTs for sums of observations from, respectively,
LRD, ND
and SRD processes over possibly non-rectangular sampling regions,
replacing $\sigma_{n}$ in (\ref{1.4}) by concise approximations which indicate the differing
rates of convergence in different situations, and highlighting an intricate
interplay between the spatial dependence structure and the geometry of the
sampling regions.

It will be observed
that the asymptotic variance of the centered sum
$[S_n -ES_n]$ shows a very complex pattern of
interactions among (i) the (effective) rate of decay of the coefficients
$\alpha(\cdot)$, (ii)~the sample size $N_n$ or equivalently,
the volume of the sampling region $R_n$,
and (iii) the shape of the
sampling region $R_n$. For simplicity of exposition, suppose for the
time being that $\alpha(\mathbf{i})$ decays as $\|\mathbf{i}\|
^{-\beta}$ as
$\|\mathbf{i}\|\to\infty$ and that the volume of $R_n$ grows
at rate $c_0\lambda_n^d$ for some $c_0\in(0,\infty)$ and $\lambda
_n\to
\infty$. Then, square-summability of the $\alpha(\cdot)$'s implies that
$\beta>d/2$ and it can be shown that
$\beta\in(\frac{d}{2},d)$ leads to the case of PSD.
In this case, we show that
%
\begin{equation}
\label{psd-lm}
\bigl(\lambda_n^{3d-2\beta} \bigr)^{-1/2}[S_n
- ES_n] \mathop{\to}^d N\bigl(0,\sigma_{\mathrm{psd}}^2
\bigr),
\end{equation}
where $\sigma_{\mathrm{psd}}^2$ depends on certain limiting
characteristics of
the co-efficients $\alpha(\cdot)$ and the shape
of the sampling region $R_n$ (cf. Theorem~\ref{th32} below).
Note that in the PSD case,
that is, for $\beta\in(d/2,d)$,
the scaling sequence is given by
${\lambda}_{n}^{(3d-2\beta)/2}$ which is of a larger
order of magnitude than $c_0^{1/2}{\lambda}_{n}^{d/2}$,
the square root of the volume of the sampling
region $R_n$. Thus, in the PSD case,
the variance of the sum grows at a rate
faster than the usual rate $N_n^{1/2}$ (since
the sample size $N_n$ here also grows at the rate
$c_0{\lambda}_{n}^{d}$). Further, the limiting
variance of the sum does \textit{not} depend on
the values of $\alpha(\mathbf{i})$ for $\mathbf{i}$ in any
given bounded neighborhood of the origin.

In the ND case, the sum shows
a very different limit behavior that critically
depends on $\beta$ as well as
on the values of the coefficients $\alpha(\mathbf{i})$,
for both small as well as large values of $\|\mathbf{i}\|$.
Indeed, for $\beta\in(d, d+1/2)$, a normal limit similar to
(\ref{psd-lm}) holds, albeit with a different asymptotic
variance, which now depends on properties of both $R_n$
and $R_n^c$. On the other hand, for $\beta$ beyond the
critical level $d+1/2$, the edge-effect of the sampling region $R_n$
becomes asymptotically dominant in dimensions $d\geq2$,
which in turn determines the asymptotic distribution.
The corresponding scaling sequence
is now given by ${\lambda}_{n}^{(d-1)/2}$
which, quite surprisingly, no longer depends on
the values of $\beta\in(d+1/2,\infty)$, that is, on the
rate of decay of the coefficients $\alpha(\mathbf{i})$. Thus, in
dimensions $d\geq2$, the slowest possible
rate for the variance of the sum in the ND case is
given by ${\lambda}_{n}^{(d-1)}$ for all $\beta>d+1/2$.
This may be contrasted with the one dimensional case,
where the edge-effect is asymptotically negligible and
the growth rate of the variance of the sum can be very slow
(e.g., $\mathrm{O}({\lambda}_{n}^{3-2\beta})$ with $\beta$ close to $d+1/2=3/2$).
See Section~\ref{sec4} for full details.
For the sake of completeness, we also prove a CLT
for the SRD case. Here the sum has the
usual rate of $N^{1/2}$ and the limiting variance
depends on the $\alpha(\mathbf{i})$ only through
their sum $A$, agreeing with a familiar result
in the time series case $d=1$
(cf. Section~\ref{sec5}). The following Table~\ref{tab1} summarizes
different limit behavior of the sum under PSD, ND and SRD
for $d\geq2$.

\begin{table}\label{tab1}
\caption{A summary of the limit behavior of the sum $S_n$
under PSD, ND and SRD when $\alpha(\mathbf{i}) = c_1\|\mathbf{i}\|
^{-\beta}$
for $\|\mathbf{i}\|>c_2$ and when $\operatorname{vol.}(R_n)\sim c_3
\lambda_n^d$,
for some constants $c_1,
c_2,c_3\in(0,\infty)$ and
for some $\lambda_n\to\infty$, where $\operatorname{vol.}(R_n)$
denotes the volume of $R_n$.
Note that here the sample size $N_n\sim\operatorname{vol.}(R_n)$
and $A=\sum_{\mathbf{i}\in\mathbb{Z}^d} \alpha(\mathbf{i})$. The
cases ND-EE and ND-NEE
in the first column
stand for the ND case with- and without-edge effects, respectively}
\begin{tabular*}{\tablewidth}{@{\extracolsep{\fill}}llll@{}}
\hline
&  &\multicolumn{2}{l@{}}{Effects on limit variance}\\[-6pt]
&   & \multicolumn{2}{l@{}}{\hrulefill} \\
&Growth rate of $\operatorname{Var}(S_n)$& Coefficients $\alpha(\mathbf{i})$ & Irregular shape of $R_n$\\
\hline
PSD & $\lambda_n^{(3d-2\beta)}, \beta\in(\frac{d}{2},d)$
& Tail behavior   at infinity
& Geometry of $R_n$\\
&&  &\\[3pt]
ND-NEE& $\lambda_n^{(3d-2\beta)}, \beta\in(d,d+\frac{1}{2})$ &
Tail behavior at  infinity
& Geometry of $R_n$ and $R_n^c$\\
&&and $A=0$ &\\[3pt]
ND-EE& $\lambda_n^{(d-1)}, \beta\in(d+\frac{1}{2},\infty)$ &
$A=0$, but not on tail  &
Geometry of $\partial R_n$, the \\
&&  behavior& boundary of $R_n$\\[3pt]
SRD & $N_n$, $\beta>d$ & Only on $A$& None\\
\hline
\end{tabular*}
\end{table}

The rest of the paper is organized as follows. In Section~\ref{sec2},
we consider the spatial linear process allowing non-identically
distributed (but homoscedastic) errors and describe an
asymptotic framework that can accommodate a large class
of sampling regions of non-standard shapes. In Section~\ref{sec2},
we also state a regular variation condition on the
coefficients that, in particular, allows
the coefficients to have different rates of decay
along different directions, and
give some examples to illustrate the
scope of the formulation.
In Sections~\ref{sec3} and \ref{sec4}, we
establish the limit distribution of the sum under PSD
and ND, respectively.
We prove the CLT in the SRD case
in Section~\ref{sec5}.
Proofs of the main results are given
in Section~\ref{sec6}. Here we also present a very general
version of the CLT for a spatial linear process
observed on bounded regions that
may be of independent interest.

\section{The theoretical framework}\label{sec2}

In Section~\ref{sec21}, we specify the spatial linear process
and in Section~\ref{sec22}, we give a formulation for the
sampling regions $R_n$. In Section~\ref{sec23}, we
introduce the regularity conditions on the coefficients
$\alpha(\mathbf{i})$ and give some illustrative examples
in Section~\ref{sec24}. Under these conditions, it is possible
to determine the exact order of the variance
term $\sigma_n^2$ and derive explicit expressions
for the asymptotic variance.

\subsection{Spatial linear processes}\label{sec21}
We define a spatial
linear process $\{Z(\cdot)\}$ as:
%
\begin{equation}\label{li-p}
Z(\mathbf{i}) = \mu+\sum_{\mathbf{j}\in\mathbb{Z}^d} \alpha (\mathbf{i}-
\mathbf{j}) \varepsilon(\mathbf{j}), \qquad\mathbf{i}\in\mathbb{Z}^d,
\end{equation}
where $ \{\varepsilon(\mathbf{j}) \dvt \mathbf{j}\in\mathbb{Z}^d\}$ is a
collection of
independent zero mean random variables with
common variance 1 (w.l.o.g.),
the $\{\varepsilon(\mathbf{j})^2 \dvt  \mathbf{j}\in\mathbb{Z}^d\}$
are uniformly integrable
and $\{\alpha(\mathbf{j})\dvt \mathbf{j}\in\mathbb{Z}^d\}$ is a
sequence of
real numbers satisfying $\sum_{\mathbf{i}\in\mathbb{Z}^d} |\alpha
(\mathbf{i})|^2<\infty$.

\subsection{Sampling regions}\label{sec22}

Next, we specify the structure of the sampling region $R_n$.
Let
$R_0$ be an open connected subset of $(-1/2,1/2]^d$ containing
the origin.
We regard $R_0$
as a ``prototype'' of the sampling region $R_n$. Let $\{\lambda_n\}$ be a
sequence of positive numbers such that $\lambda_n
\to\infty$
as $n\rightarrow\infty$.
We assume that the sampling region $R_n$ is
obtained by ``inflating'' the set $R_0$ by the scaling
factor $\lambda_n$ (cf. Lahiri \textit{et al.} \cite{Lahetal99}),
that is,
%
\begin{equation}\label{eq22}
R_n=\lambda_nR_0.
\end{equation}
Since the origin is assumed to lie in $R_0$, the shape of $R_n$ remains
the same for different values of $n$.
To avoid pathological cases, we assume that the boundary
$\partial R_0$ of $R_0$ has $d$-dimensional Lebesgue measure
zero. A stronger version of this condition will be needed
for the ND case, which is stated as condition (C.3) in Section~\ref{sec23} below.
The (stronger) boundary condition holds for
most regions $R_n$ of practical
interest, including common convex subsets of $\mathbb{R}^d$,
such as spheres, ellipsoids, polyhedrons, as well as for many
non-convex star-shaped sets in $\mathbb{R}^d$. (Recall that a set
$A\subset\mathbb{R}^d$ is called {\it star-shaped} if for any $x\in
A$, the
line segment joining $x$ to the origin lies in $A$.)
The latter class of sets may have
fairly irregular shapes (cf. Sherman and Carlstein \cite{SheCar94}
and Lahiri \cite{Lah99}). Some practical applications and studies
involving sampling regions that satisfy the regularity
conditions above are given by
the wheat yield data of Mercer and Hall \cite{MerHal11} on
agricultural field trials, the coal ash data of Gomez and Hazen \cite{GomHaz70} from Mining, and the cancer mortality counts data
of Riggan \textit{et al.} \cite{Rigetal87} from Epidemiology, among others.

\subsection{Regularity conditions}\label{sec23}
For $\mathbf{x}=(x_1,\ldots,x_d)'\in{\mathbb{R}}^d$,
let $\|\mathbf{x}\|= (x_1^2+\cdots+x_d^2)^{1/2}$
and let $\lfloor{\mathbf{x}}\rfloor = (\lfloor{x_1}\rfloor,\ldots,\lfloor{x_d}\rfloor)'$
where
$\lfloor{y}\rfloor$ denotes the integer part of a real number $y$.
For $\delta\in(0,\infty)$, $\mathbf{x}\in\mathbb{R}^d$ and
$A\subset\mathbb{R}^d$, let
$A^\delta= \{\mathbf{y}\in\mathbb{R}^d \dvt  \|\mathbf{z}-\mathbf
{y}\| \leq\delta$
for some $\mathbf{z}\in A\}$ and $d(\mathbf{x},A) \equiv d(A,\mathbf{x})
=\min\{\|\mathbf{z}-\mathbf{x}\|\dvt  \mathbf{z}\in A\}$.
Let
$
\gamma(t) =\max\{ | \alpha(\lfloor{\mathbf{u}t}\rfloor)| \dvt \|
\mathbf{u}\|=1\}, t>0$.
With this notation, we are now ready to state the regularity conditions.

\begin{enumerate}[(C.3)]
\item[(C.1)]
Suppose that
\[
\gamma(t) = t^{-\beta} L(t),\qquad  t>0,
\]
for some $\beta>d/2$
and some function $L\dvtx  (0,\infty) \to[0,\infty)$ that is slowly varying
at infinity in the sense that $L(\cdot)$ is bounded on any bounded
subinterval of
$(0,\infty) $ and $\lim_{t \to\infty} \sup\{L(at)/L(t) \dvt  a\in
[a_0,a_1]\} =1$
for any $0<a_0<a_1<\infty$ (cf. Taqqu \cite{Taq74}).
\item[(C.2)]
Let $g_t(\mathbf{x}) = \alpha(\lfloor{t \mathbf{x}}\rfloor)/\gamma
(t)$, $\mathbf{x}\in\mathbb{R}^d$, $t>0$.
Suppose that there exists a function ${g}_{\infty}\dvtx \mathbb{R}^d\to
\mathbb{R}$ such
that for every $\delta\in(0,\infty)$,
\[
\int_{\{\mathbf{x}\in\mathbb{R}^d\dvt \|\mathbf{x}\|\geq\delta\}} \bigl| g_t(\mathbf{x}) -
{g}_{\infty}(\mathbf{x}) \bigr|^{b}\,\mathrm{d}\mathbf{x}\to0  \qquad\mbox{as } t\to
\infty,
\]
where $b = b(d,\beta) = 2$ if $\beta\leq d$ and $b=1$ otherwise.
%
\item[(C.3)]
For any measurable function $f\dvtx  [0,\infty) \to[0,\infty)$, there exists
$C_f\in(1,\infty)$ such that
\[
\int_{[\partial R_0]^\varepsilon} f \bigl(d(\mathbf{x},\partial R_0)
\bigr) \nu(\mathrm{d}\mathbf{x}) \leq C_f \int_0^\varepsilon
f(t) \,\mathrm{d}t  \qquad \mbox{for all } 0<\varepsilon<C_f^{-1},
\]
where
$\nu$ is the Lebesgue measure on $\mathbb{R}^d$.
\end{enumerate}

Condition (C.1) requires that the radial maximum of the
collection of coefficients $\{\alpha(\mathbf{i})\dvt  \mathbf{i}\in
\mathbb{Z}^d\}$
be regularly varying (at infinity).
The requirement that $\beta>d/2$ in (C.1) is imposed to ensure
that $\sum_{\mathbf{i}\in\mathbb{Z}^d}\alpha(\mathbf{i})^2
<\infty$.
Condition (C.2) is a weak form of spatial regular variation
condition on the $\alpha(\mathbf{i})$.
It is weaker than assuming directional
separability of the coefficients, and allows for differential
rates of decay along different directions. See the examples below.
It is a variant of the standard form of regular variation
that requires the function $g_t(\cdot)$ to satisfy
(cf. Section~5.4, Resnick \cite{Res87})
%
\begin{equation}\label{reg-v}
\lim_{t \to\infty} g_t(\mathbf{x}) = \|\mathbf{x}
\|^{-\beta}a\bigl(\mathbf{x}/\|\mathbf{x}\|\bigr)  \qquad\mbox{for all } \mathbf{x}\in
\mathbb{R}^d, \mathbf{x}\neq\mathbf{0},
\end{equation}
for some function $a(\cdot)$ on the unit disc
$\{\mathbf{x}\in\mathbb{R}^d\dvt  \|\mathbf{x}\|=1\}$. In this case,
the limit function ${g}_{\infty}$ is given by
\[
{g}_{\infty}(\mathbf{x})= \|\mathbf{x}\|^{-\beta}a\bigl(\mathbf{x}/\|
\mathbf{x}\|\bigr)\mathbh{1}(\mathbf{x}\neq\mathbf{0}),\qquad \mathbf{x}\in
\mathbb{R}^d,
\]
where $\mathbh{1}(\cdot)$ denotes the
indicator function. By comparison, condition (C.2) requires
convergence of $g_t$ to ${g}_{\infty}$ in $L^b$.
Conditions (C.1) and (C.2) together quantify the behavior of the
function $\alpha(\mathbf{i})$ for large $\|\mathbf{i}\|$ which
plays an important role in determining the
form of the asymptotic variance of the sum under PSD and ND.

Condition (C.3) is a regularity condition on the boundary of
the prototype set $R_0$
which is equivalent to requiring that
\[
\int_{[\partial R_0]^\varepsilon} f \bigl(d(\mathbf{x},\partial R_0)
\bigr) \nu(\mathrm{d}\mathbf{x}) =\mathrm{O} \biggl( \int_0^\varepsilon
f(t)\, \mathrm{d}t \biggr)\qquad\mbox{as } \varepsilon\downarrow0
\]
for each non-negative measurable $f$. We need this condition to hold
for $f\equiv1$ and for $f(t) = t^{-b}L^2(t)$ for certain values of
$b=b(\beta)$ (cf. the proof of Theorem~\ref{th41} below). In particular, when
$f\equiv1$, this reduces to the condition
%
\begin{equation}\label{bd-cd}
\nu \bigl((\partial R_0)^\varepsilon \bigr) = \mathrm{O}(\varepsilon)
\qquad \mbox{as } \varepsilon\downarrow 0,
\end{equation}
which is satisfied by most sampling regions
of common interest (cf. Section~\ref{sec22}). For $d=2$, a~sufficient condition
is
that the boundary of $R_0$ is delineated by a rectifiable curve
of a finite length.
We shall use (C.3) for proving the
results \textit{only} in the ND case where more precise
information on the bounadry is needed to determine the
asymptotic variance.

Next, we give a few examples to illustrate
the range of spatial dependence
covered by the regularity conditions above.

\subsection{Examples}\label{sec24}

\begin{exa}[(Isotropic spatial linear processes)]\label{ex21}
Let
\[
\alpha(\mathbf{i}) = a\bigl(\|\mathbf{i}\|\bigr)\bigl(1+\|\mathbf{i}\|\bigr)^{-\beta},
\qquad \mathbf{i}\in\mathbb{Z}^d,
\]
for some bounded function $a\dvtx  [0,\infty) \to
\mathbb{R}$, where $\beta\in(\frac{d}{2},\infty)$. Also, suppose that
$a(t)\to c_0\neq0$ as $t\to\infty$. Then, using the fact that
%
\begin{equation}\label{ip-ineq}
\sup \bigl\{ \bigl|\bigl\|\lfloor{t\mathbf{u}}\rfloor\bigr\| - t\|\mathbf{u}\| \bigr| \dvt  \mathbf{u}
\in\mathbb{R}^d\setminus\{\mathbf{0}\} \bigr\} \leq\sqrt{d},
\end{equation}
it is easy to see that condition (C.1) holds with $\gamma(t)
= t^{-\beta}L(t)$
where $L(t)\to|c_0|$ as $t\to\infty$. Further, using (\ref{ip-ineq}),
one can show that for any $\eta>0$,
\[
g_t(\mathbf{x}) \to\frac{c_0}{|c_0|} \|\mathbf{x}
\|^{-\beta} \equiv{g}_{\infty}(\mathbf{x})\qquad \mbox{as }t \to\infty \mbox{ for all }\|\mathbf{x}\|>\eta
\]
and
\[
\bigl|g_t(\mathbf{x}) - {g}_{\infty}(\mathbf{x}) \bigr| \leq
C_1 \| \mathbf{x}\|^{-\beta} \qquad\mbox{for all } t> C_1,
\]
for some constant $C_1 \equiv C_1(\eta, \beta)\in(0,\infty)$.
Hence, condition (C.2) holds.

This gives an example of an
``isotropic'' spatial linear process where the coefficients
$\alpha(\mathbf{i})$ have
an identical rate of decay
in all directions.
\end{exa}

\begin{exa}[(A class of anisotropic spatial linear
processes)]\label{ex22}
Suppose that $O$ is a $d\times d$
orthonormal matrix with rows $\mathbf{o}_i'$, $i=1,\ldots,d$.
Let $\phi_i(\mathbf{x}) = |\mathbf{o}_i'\mathbf{x}|/\|\mathbf{x}\|$,
$\mathbf{x}\in\mathbb{R}^d\setminus\{\mathbf{0}\}$.
Let $\delta\in(0, \frac{1}{\sqrt{d}})$.
Suppose that
\[
\alpha(\mathbf{x}) = \prod_{i=1}^d \bigl\{\bigl|\mathbf{o}_i' \mathbf{x}\bigr|^{-a_i} \mathbh{1}
\bigl(\phi_i(\mathbf{x})> \delta \bigr) \bigr\}\qquad \mbox{for all }
\mathbf{x} \in\mathbb{R}^d\setminus\Gamma
\]
for some $a_1,\ldots,a_d \in[0,\infty)$ with $a_1+\cdots+a_d >d/2$
and for some open neighborhood $\Gamma$ of the origin. There
is no
restriction on the definition of the function
$\alpha(\mathbf{x})$ on $\Gamma$.

It is easy to check that
\begin{eqnarray*}
\gamma(t) &=& \sup \bigl\{ \bigl|\alpha\bigl(\lfloor{\mathbf{u}t}\rfloor \bigr) \bigr| \dvt
\|\mathbf{u}\| =1 \bigr\}
\\
&=& \sup \bigl\{ \bigl|\alpha\bigl(\bigl\lfloor{O'\mathbf{u}t}\bigr
\rfloor\bigr) \bigr| \dvt  \|\mathbf{u}\| =1 \bigr\}
\\
&=& \sup \Biggl\{\prod_{i=1}^d \bigl|
\lfloor{u_i t}\rfloor \bigr|^{-a_i} \mathbh{1} \bigl(
\phi_i\bigl(\lfloor{t\mathbf{u}}\rfloor\bigr)> \delta \bigr) \dvt  \|
\mathbf{u}\| =1 \Biggr\}\bigl(1+\mathrm{o}(1)\bigr)
\\
&=& c_0 t^{-(a_1+\cdots+a_d)} \bigl(1+\mathrm{o}(1)\bigr)\qquad\mbox{as }t \to\infty,
\end{eqnarray*}
for some
$c_0\in(0,\infty)$.

Next, let $D= \{ \mathbf{0}\}\cup\{\mathbf{y}\dvt  \phi_i(\mathbf{y})
=\delta$ for
$i=1,\ldots,d\}$. Note that the $d$-dimensional Lebesgue
measure of $D$ is zero. And, for any $\mathbf{x}\notin D$, 
\begin{eqnarray*}
g_t(\mathbf{x}) &=& \alpha\bigl(\lfloor{t\mathbf{x}}\rfloor\bigr)/\gamma (t)
\\
&=& \gamma(t)^{-1} \prod_{i=1}^d
\bigl|\bigl\lfloor{t \mathbf{o}_i'\mathbf{x}}\bigr\rfloor\bigr|^{-a_i} \mathbh{1} \bigl(\phi_i\bigl(\lfloor{t\mathbf{x}}
\rfloor\bigr)>\delta \bigr)
\\
&=& \Biggl[ c_0^{-1} \prod_{i=1}^d
\bigl\{\bigl|\mathbf{o}_i'\mathbf{x}\bigr|^{-a_i}
\mathbh{1} \bigl(\phi_i(\mathbf{x})> \delta \bigr) \bigr\} \Biggr]
\bigl(1+\mathrm{o}(1)\bigr)\qquad\mbox{as }t \to\infty,
\\
&\equiv& {g}_{\infty}(\mathbf{x}) \bigl(1+\mathrm{o}(1)\bigr)\qquad\mbox{as }t \to
\infty.
\end{eqnarray*}
Thus, the point-wise limit of the functions $g_t(\cdot)$ exists
for all $\mathbf{x}\notin D$. Now using the Dominated Convergence
Theorem (DCT),
it is easy to check that for any $\eta>0$,
\[
\lim_{t \to\infty} \int_{\{\|\mathbf{x}\|\geq\eta\}} \bigl|g_t(
\mathbf{x}) - {g}_{\infty}(\mathbf{x}) \bigr|^b \,\mathrm{d}\mathbf{x}=0.
\]
Thus, conditions (C.1) and (C.2) are satisfied by the
coefficients generated by the
function $\alpha(\cdot)$.

Note that in this example,
the coefficients
$\alpha(\mathbf{i})$ are {\it zero} whenever $|\phi_i(\mathbf{i})|
\leq\delta$
for some $i\in\{1,\ldots,d\}$
but they are non-zero for all $\mathbf{x}$
such that $\phi_i(\mathbf{x}) >\delta$ for all $i=1,\ldots,d$.
Since $\delta$ is small, the latter condition is satisfied
in a conic region in each of the $2^d$ quadrants.
Thus, this gives an example of an anisotropic
spatial process.
The maximal rate of decay of $\alpha(\cdot)$ over the set
$\{\mathbf{i}\in\mathbb{Z}^d\dvt  |\phi_i(\mathbf{i})| > \delta$ for all
$i=1,\ldots, d\}$ can vary with the choice of $a_1,\ldots, a_d$,
allowing all possible types of long-range (as well as
short-range) dependence.
Also, note that here the ND case can be realized by
a suitable choice of $\alpha(\mathbf{i})$ for $\mathbf{i}\in\Gamma$.
The rate of convergence of the sum in this example depends
only on a {\it single} parameter, namely,
the combined exponent $\beta= a_1+\cdots+a_d$;
Individual $a_i$'s do not have an impact.
\end{exa}

%

\begin{exa}[(Spatial linear processes with non-uniform directional decay
rates)]\label{ex23}
Let $\mathcal{I}$ be a~finite set and
let $\{\mathbf{o}_{i}\dvt  i\in\mathcal{I}\}\subset\{\mathbf{x}\in
\mathbb{R}^d\dvt
\|\mathbf{x}\|=1\}$. Here we suppose that the $\mathbf{o}_i$'s are distinct
but they are not necessarily orthogonal.
Let $\phi_i(\cdot)$ be
as in Example~\ref{ex22}, that is,
$\phi_i(\mathbf{x}) = |\mathbf{o}_i'\mathbf{x}|/\|\mathbf{x}\|$,
$\mathbf{x}\in\mathbb{R}^d\setminus\{\mathbf{0}\}$, $i\in\mathcal{I}$.
Define $\psi_i(\mathbf{x}) = \phi_i(\mathbf{x}) \mathbh{1}(\phi
_i(\mathbf{x}) >\delta_i)$
for some $\delta_i\in(0,1)$, $i\in\mathcal{I}$. Suppose that
\[
\alpha(\mathbf{x}) = \sum_{i\in\mathcal{I}} \frac{\psi
_i(\mathbf{x})}{
1+\|\mathbf{x}\|^{a_{i}} }
\qquad \mbox{for all } \mathbf{x} \in\mathbb{R}^d\setminus\Gamma,
\]
where
$\{a_{i} \dvt   i\in\mathcal{I}\}\subset(0,\infty)$ and where
$\Gamma$ is an open neighborhood of the origin, as in Example~\ref{ex22}.
Let $a_0=
\min\{a_{i}\dvt  i\in\mathcal{I}\}$ and let
$\mathcal{I}_0 =\{i\in\mathcal{I}\dvt  a_{i} = a_0\}$.
Note
that for any $i\in\mathcal{I}$ and any
$\mathbf{x}\in\mathbb{R}^d\setminus\{\mathbf{0}\}$,
\[
\phi_i\bigl(\lfloor{t\mathbf{x}}\rfloor\bigr) =
\phi_i(\mathbf {x}) \bigl(1+\mathrm{o}(1)\bigr)\qquad \mbox{as }t \to\infty,
\mbox{provided }\phi_i(\mathbf{x}) >0,
\]
and
\[
\bigl\|\lfloor{t\mathbf{x}}\rfloor\bigr\| =
t\|\mathbf{x}\| \bigl(1+\mathrm{o}(1)\bigr)\qquad
\mbox{as }t \to\infty.
\]
Hence, it follows that
\[
\gamma(t) = c_1 t^{-a_0}\bigl(1+\mathrm{o}(1)\bigr) \qquad\mbox{as }t \to
\infty,
\]
where $c_1 = \sup\{\sum_{i\in\mathcal{I}_0} \psi_i(\mathbf{u}) \dvt
\|\mathbf{u}\|=1\}$.

Next, we identify the limit function ${g}_{\infty}(\cdot)$.
By arguments as above, it follows that for any $\mathbf{x}\neq\mathbf{0}$
with $\psi_i(\mathbf{x}) > 0$ for some $i\in\mathcal{I}_0$,
\begin{eqnarray*}
g_t(\mathbf{x}) &=& \gamma(t)^{-1} \biggl[ \sum
_{i\in\mathcal{I}} \frac{ \psi
_i(\lfloor{t\mathbf{x}}\rfloor)}{
1+ \|\lfloor{t\mathbf{x}}\rfloor\|^{a_{i}}} \biggr]
\\
&=&\frac{\sum_{i\in\mathcal{I}_0} \psi_i(\mathbf{x})}{c_1 \|
\mathbf{x}\|^{a_0}} \bigl(1+\mathrm{o}(1)\bigr) \qquad\mbox{as }t \to\infty, \mbox{a.e.}
\end{eqnarray*}
Since for any $i$, the set $\{\mathbf{x}\in\mathbb{R}^d\setminus\{
\mathbf{0}\} \dvt
\phi_i(\mathbf{x}) = \delta_i\}$ has $d$-dimensional Lebesgue
measure zero, it follows that
$g_t \to {g}_{\infty}$ {as} $t \to\infty$ (a.e.),
where
\[
{g}_{\infty}(\mathbf{x}) = \frac{\sum_{i\in\mathcal{I}_0}
\phi_i(\mathbf{x})}{c_1 \|\mathbf{x}\|^{a_0}}, \qquad \mathbf{x}\in
\mathbb{R}^d\setminus\{\mathbf{0}\}.
\]
Now using the DCT, one can show that for any $\eta>0$,
\[
\lim_{t \to\infty} \int_{\{\|\mathbf{x}\|\geq\eta\}}  \bigl|g_t(
\mathbf{x}) - {g}_{\infty}(\mathbf{x}) \bigr|^b \,\mathrm{d}\mathbf{x}=0.
\]
Thus, conditions (C.1) and (C.2) are satisfied by the
coefficients generated by the
function
$\alpha(\cdot)$.

Note that for $\delta_i$ close to 1, the $i$th
component $\frac{\psi_i(\mathbf{x})}{1+\|\mathbf{x}\|^{a_{i}}}$
in $\alpha(\mathbf{x})$ takes non-zero values\vspace*{2pt} in a thin
cone around $\mathbf{o}_i$ and it may or may not
intersect with the $j$th cone (for any given $j\neq i$)
depending on the relative magnitudes of $\delta_i$ and $\delta_j$
and the angle between $\mathbf{o}_i$ and $\mathbf{o}_j$. As a result,
with different choices of $\mathbf{o}_i, \delta_i$ and $a_i$
for $i\in\mathcal{I}$, the coefficients
$\alpha(\mathbf{i})$ here may have different rates of decay along the
directions $\mathbf{o}_i$ for $i\in\mathcal{I}\setminus\mathcal
{I}_0$, allowing
{\it any combinations of
short- and long-range dependent rates along different directions}.
However, the limit distribution
of the sum depends only on $a_0$
which is the minimum of the exponents $\{a_i, i\in\mathcal{I}\}$.
\end{exa}



In the next section, we describe the limit behavior of the
sum $S_n$ depending on the rate of decay $\beta$ in (C.1).



\section{Results under PSD}\label{sec3}
From the proofs given in Section~\ref{sec6},
it follows that for $\beta\in(d/2, d)$,
the variance of the sum $S_n$
grows at the rate ${\lambda}_{n}^{3d-2\beta}L({\lambda}_{n})^2$
and hence, the correct scaling factor sequence is given by
${\lambda}_{n}^{({3d-2\beta})/{2}}L({\lambda}_{n})$. Since $\beta
\in(d/2, d)$,
this scaling sequence grows at a rate {\it faster} than
square-root of the sample size $|N_n|^{1/2}\sim[{\lambda}_{n}^{d}
\operatorname{vol.}(R_0)]^{1/2}$, where $\operatorname{vol.}(B)$
denotes the volume (i.e.,
the Lebesgue measure) of a Borel set $B$ in $\mathbb{R}^d$.
As a result, for $\beta\in(d/2, d)$,
$N_n^{-1}\sigma_n^2\to\infty$
and the spatial process $\{Z(\cdot)\}$ exhibits PSD.

To describe the limit distribution of the centered and scaled
sum under PSD, define
%
\begin{equation}\label{Gi}
{G}_{\infty}(\mathbf{x}) = 
\int_{R_0}
{g}_{\infty}(
\mathbf{y}-\mathbf{x}) \,\mathrm{d}\mathbf{y}, \qquad\mathbf{x}\in\mathbb{R}^d.
\end{equation}
%

The first result
shows that the function ${G}_{\infty}$ is
well defined on all of $\mathbb{R}^d$ for $\beta\in(0,d)$.

\begin{prop}\label{pr3.1}
Suppose that conditions \textup{(C.1)} and \textup{(C.2)}
hold for some $\beta\in(0,d)$. Then
the integral in (\ref{Gi}) exists and is finite for all
$\mathbf{x}\in\mathbb{R}^d$.
\end{prop}

The function ${G}_{\infty}$ determines the asymptotic variance of the
sum $S_n$ for $\beta\in(d/2, d)$.
We make this precise in the
following result that gives the limit distribution of
$S_n$ under PSD.

\begin{theorem}\label{th32}
Let $\{Z(\cdot)\}$ be the linear process given by
(\ref{li-p}) such that conditions \textup{(C.1)} and \textup{(C.2)} hold for some $\beta
\in(d/2,d)$. Then
${G}_{\infty}\in L^2(\mathbb{R}^d)$ and
%
\begin{equation} \label{clt-psd}
\frac{[S_n-ES_n]}{
 [{\lambda}_{n}^{{3d-2\beta}} L({\lambda}_{n})^2 ]^{1/2}}  \mathop{\to}^d  N \biggl(0, \int_{\mathbb{R}^d}
{G}_{\infty}^2(\mathbf {x}) \,\mathrm{d}\mathbf{x} \biggr)\qquad \mbox{as }n\to
\infty.
\end{equation}
\end{theorem}

Theorem~\ref{th32} shows that for $\beta\in(d/2, d)$, the growth
rate of (the variance of) the sum $S_n$ is
$[{\lambda}_{n}^{{3d-2\beta}} L({\lambda}_{n})^2]^{1/2}$,
which is of a larger order than the usual order
$N_n^{1/2}$ .
Since ${G}_{\infty}$ is continuous, the asymptotic
variance is non-zero if
\[
{G}_{\infty}(\mathbf{x}_0) \neq0 \qquad\mbox{for some }
\mathbf{x}_0\in \mathbb{R}^d.
\]
From (\ref{clt-psd}), also note that the limiting
variance of $S_n$ depends on the
prototype set $R_0$ as well as the function
${g}_{\infty}$ of condition (C.2). Thus, unlike the time series case,
the geometry of the sampling region plays an important role in
the spatial case under PSD.

\section{Results under ND}\label{sec4}
When $Z(\mathbf{i})$ is ND,
%
it can be shown (cf. the proofs of Theorems~\ref{th41}, \ref{th43} and \ref{thm44})
that
$N_n^{-1}\sigma_n^2\to0$.
The limit behavior of the
sum $S_n$ in the ND case critically depends on the
behavior of the terms
\[
\theta_n(\mathbf{i}) = \sum_{\mathbf{j}\in[R_n-\mathbf{i}]\cap
\mathbb{Z}^d}
\alpha(\mathbf{j}), \qquad \mathbf{i}\in\mathbb{Z}^d,
\]
in a shrinking neighborhood of the set, $\partial R_n$,
the boundary of $R_n$. In the parlance of spatial statistics,
this represents an instance of {\it edge effect} (cf. Cressie \cite{Cre93})
that may have a non-trivial effect on the limit behavior
of the sum. Indeed,
depending on the relative orders of contributions
from the boundary terms and the non-boundary terms, we
get different growth rates for the sum
in the ND case. Further, the limiting variances are also
different. For clarity of exposition, we present the two
subcases of the ND case separately.

\subsection{ND with asymptotically negligible edge effects}\label{sec41}
First, we consider the relatively simple case where
the contribution from the boundary $\theta_n(\mathbf{i})$'s
is asymptotically negligible. Suppose that $\beta>d$,
so that $\sum_{\mathbf{i}\in\mathbb{Z}^d} |\alpha(\mathbf{i})|
<\infty$,
and that $A\equiv \sum_{\mathbf{i}\in\mathbb{Z}^d} \alpha(\mathbf
{i}) =0$.
In this case, it will be shown that the
asymptotic distribution of the sum depends on
the $\alpha(\mathbf{i})$ only through (an analog of)
the function ${G}_{\infty}$ of (\ref{Gi}). However, Proposition~\ref{pr3.1} no longer holds, that is, the function ${G}_{\infty}$ of (\ref{Gi})
may not be well defined for all $\mathbf{x}\in\mathbb{R}^d$
for $\beta>d$.
To appreciate why, consider the special case
where ${g}_{\infty}(\mathbf{x}) = \|\mathbf{x}\|^{-\beta}\mathbh
{1}(\mathbf{x}\neq\mathbf{0})$
(cf. (\ref{reg-v})). In this case,
$\int_{\{\delta\leq\|\mathbf{y}\|\leq1\}}{g}_{\infty}(\mathbf{y})
=\mathrm{O}(\delta^{d-\beta}) = \mathrm{o}(1)$ as $\delta\downarrow0$ for all $\beta
\in(d/2,d)$, but the
integral blows up for $\beta> d$. As a result, the limit
in (\ref{Gi}) may not exist {\it for all}
$\mathbf{x}\in\mathbb{R}^d\setminus\{\mathbf{0}\}$ in the case
$\beta>d$.
However, using the condition $A=0$, we can define
${G}_{\infty}(\mathbf{x})$ (and a suitable variant of it)
for a restricted set of $\mathbf{x}$'s that
would be adequate for our purpose.

To that end, note that the set $[R_0\cup\partial R_0]^c$ is
open and hence, for all $\mathbf{x}\in[R_0\cup\partial R_0]^c$,
there exists a $\eta=\eta(\mathbf{x})>0$\vspace*{-3pt}
such that
\[
B(\mathbf{x};\eta)\subset[R_0\cup\partial R_0]^c,
\]
where
$B(\mathbf{x};\eta) \equiv \{\mathbf{y}\in\mathbb{R}^d\dvt  \|\mathbf
{y}-\mathbf{x}\| <\eta\}$ denotes
the open ball of radius $\eta$ around $\mathbf{x}$.
As a consequence,
for $\mathbf{x}\in[R_0\cup\partial R_0]^c$,
$B(\mathbf{0};\eta)\cap[R_0-\mathbf{x}]=\varnothing$\vspace*{-3pt} and
\[
\int_{R_0}\bigl|{g}_{\infty}(\mathbf{y}-\mathbf{x})\bigr|\, \mathrm{d}
\mathbf{y} = \int_{R_0-\mathbf{x}}\bigl|{g}_{\infty}(\mathbf{y})\bigr|\, \mathrm{d}
\mathbf{y} \leq \biggl[\int_{\{\|\mathbf{y}\|\geq\eta\}} \bigl|{g}_{\infty}(\mathbf
{y})\bigr|^2\, \mathrm{d}\mathbf{y} \biggr]^{1/2} \bigl[
\operatorname{vol.}(R_0) \bigr]^{1/2}<\infty
\]
whenever condition (C.2) holds with $b=2$. But,
for $\beta>d$, $b=1$ in condition (C.2). Nonetheless,
the square integrability of ${g}_{\infty}(\cdot)$
on sets of the form $B_{\eta}\equiv\{\mathbf{y}\in\mathbb{R}^d\dvt  \|
\mathbf{y}\|\geq\eta\}$,
$\eta>0$,
follows from (C.2) and the fact that
$|{g}_{\infty}(\cdot)| \leq C(\eta)$ a.e. (w.r.t. the
Lebesgue measure on $\mathbb{R}^d$) on $B_{\eta}$, for some $C(\eta
)\in(0,\infty)$
(see (\ref{gi-bd}) below).\vspace*{-3pt}
Hence,
\[
{G}_{\infty}(\mathbf{x}) 
= \int_{R_0}
{g}_{\infty}(\mathbf{y}-\mathbf{x}) \,\mathrm{d}\mathbf{y}\in \mathbb{R}
\]
for all $\mathbf{x} \in [R_0\cup\partial R_0]^c$.
By similar arguments, 
the integral $\int_{[R_0\cup\partial R_0]^c} {g}_{\infty}(\mathbf
{y}-\mathbf{x}) \,\mathrm{d}\mathbf{y}$
is well defined for all $\mathbf{x}\in R_0$.
Since $\partial R_0$
has $d$-dimensional Lebesgue measure zero, the value of the
integrals remains unchanged (with any measurable extension of
${g}_{\infty}(\cdot)$)
if $ [R_0\cup\partial R_0]^c$ is replaced by $R_0^c$
(cf. Billingsley \cite{Bil95} or
Athreya and Lahiri \cite{AthLah06}, page 49).
With this convention, define the function
${G^\dagger_{\infty}}(\cdot)$\vspace*{-6pt} as
%
\begin{equation}\label{cgi}
{G^\dagger_{\infty}}(\mathbf{x}) =\cases{
\ds\int_{R_0} {g}_{\infty}(\mathbf{y}-\mathbf{x})\, \mathrm{d}
\mathbf{y}, &\quad$\mbox{if }\mathbf{x}\in[R_0\cup\partial
R_0]^c$,
\vspace*{2pt}\cr
\ds\int_{R_0^c} {g}_{\infty}(\mathbf{y}-\mathbf{x}) \,\mathrm{d}
\mathbf{y}, &\quad$\mbox{if } \mathbf{x}\in R_0$,
\vspace*{2pt}\cr
0,  & \quad$\mbox{if }\mathbf{x}\in\partial R_0$.}
\end{equation}

For $\beta\in(d, d+1/2)$, the asymptotic distribution of the
sum depends only on the function ${G^\dagger_{\infty}}(\cdot)$, as
shown by the
following\vspace*{-3pt} result.

\begin{theorem}\label{th41}
Let $\{Z(\cdot)\}$ be the linear process given by
(\ref{li-p}) such that conditions \textup{(C.1)}--\textup{(C.3)} hold
with
$\beta\in(d,d+1/2)$ in \textup{(C.1)}. Also suppose that $A= 0$.
Then, ${G^\dagger_{\infty}}\in L^2(\mathbb{R}^d)$\vspace*{-3pt} and
%
\begin{equation}\label{clt-nd}
\frac{[S_n- ES_n]}{
 [{\lambda}_{n}^{{3d-2\beta}} L^2({\lambda}_{n}) ]^{1/2}}  \mathop{\to}^d  N \biggl(0, \int_{\mathbb{R}^d}
\bigl[{G^\dagger_{\infty}}(\mathbf{x})\bigr]^2\, \mathrm{d}
\mathbf{x} \biggr)\qquad \mbox{as }n\to\infty,
\end{equation}\vspace*{-3pt}%
where the function ${G^\dagger_{\infty}}$ is as defined in (\ref{cgi}).
\end{theorem}

Theorem~\ref{th41} shows that the asymptotic variance of the centered and
scaled sum depends on the coefficients $\alpha(\mathbf{i})$ only through
the integral of the function ${G^\dagger_{\infty}}(\cdot)^2$
over $\mathbb{R}^d$. Thus, the behavior of the $\alpha(\mathbf{i})$
for large values of $\|\mathbf{i}\|$ determines the asymptotic
variance. The exact values of the $\alpha(\mathbf{i})$
for small values of $\|\mathbf{i}\|$ have no direct effect
except for the condition $A=0$. Further, the growth rate of
the sum under the ND case is $[{\lambda}_{n}^{{3d-2\beta}} L({\lambda
}_{n})^2]^{1/2}
= \mathrm{o}(N_n^{1/2})$, which is
slower than the PSD rate and, also slower than the SRD rate, given
by $N_n^{1/2}$ (cf. Theorem~\ref{th51} below).
To compare the asymptotic variances under the ND case without edge effects
and the PSD case, note that the integrals
of ${G^\dagger_{\infty}}(\cdot)$ and ${G}_{\infty}(\cdot)$ over
$R_0^c$ are the
same and hence, the difference in the asymptotic variances
in the ND and the PSD cases comes
from the integrals of the respective functions over $R_0$.


\subsection{ND with asymptotically non-negligible edge effects}\label{sec42}

Next, consider the case where $\beta\geq d+1/2$. In this case,
we may write the variance of the sum as the sum of two
terms, one involving the sum of $\theta_n(\mathbf{i})^2\sigma^2$ for
$\mathbf{i}$
near the boundary of the sampling region $R_n$ and the other
over the rest of the $\theta_n(\mathbf{i})^2\sigma^2$. It can be shown
that the growth rate of the second term is of the order
${\lambda}_{n}^{{3d-2\beta}} L({\lambda}_{n})^2$. On the other hand, under
condition (C.3) (cf. (\ref{bd-cd})),
for any sequence $\{t_n\}\subset(0,\infty)$
with $t_n^{-1}+\lambda_n^{-1}t_n=\mathrm{o}(1)$ as
$n\to\infty$, the volume
of the $t_n$-enlargement of the boundary of $R_n$
is of the order of $\lambda_n^{d-1} t_n$. It is easy to check that
for $\beta\geq d+1/2$, this boundary term
can be of a {\it larger} order of magnitude than
${\lambda}_{n}^{{3d-2\beta}} L({\lambda}_{n})^2$. As a result,
the contribution from $\theta_n(\mathbf{i})^2$ for $\mathbf{i}$
near the boundary of the sampling region $R_n$
may become dominant
and
additional care must be taken to determine
the exact
growth rate of $\sigma_n^2$. The following example serves
to illustrate such dominating
``edge effects'' in the ND case:

\begin{exa}\label{ex42}
Suppose that $d=2$,
$R_0=(-\frac{1}{2},\frac{1}{2})\times
(-\frac{1}{2},\frac{1}{2})$ and let
\[
\alpha(i,j) = \cases{b(i)b(j), & \quad$\mbox{if }i,j\in
\mathbb{Z}\setminus\{0\}$,
\vspace*{2pt}\cr
0,  & \quad$\mbox{if }ij=0, (i,j)\neq(0,0)$,
\vspace*{2pt}\cr
-4B^2 & \quad$\mbox{if }(i,j)= (0,0)$,}
\]
where $\{b(i) \dvt  i\geq1\} \subset(0,\infty)$,
$b(-i) = b(i)$ for $i\geq1$ and $B\equiv\sum_{i=1}^\infty
b(i) \in(0,\infty)$. Further, suppose that $b(i)\sim
c_0 i^{ - \beta}$ as $i\to\infty$, for some $\beta>d+1/2=2.5$.
Then, $A=\sum_{\mathbf{i}\in\mathbb{Z}^2} \alpha(\mathbf{i}) =0$
and $\gamma(t) \sim c_0^2t^{-\beta}$ as $t \to\infty$.
Further, we may write $\sigma_n^2$ as
\[
\sigma_n^2 = \sum_{k=1}^3
\sum_{\mathbf{i}\in I_{kn}} \theta _n(
\mathbf{i})^2,
\]
where $I_{1n}=[-\lambda_n/2 + c_n, \lambda_n/2 - c_n]^2\cap\mathbb{Z}^2$,
$I_{2n}=\mathbb{Z}^2 \setminus[-\lambda_n/2 - c_n, \lambda_n/2 + c_n]^2$,
and $I_{3n} = \mathbb{Z}^2 \setminus[I_{1n}\cup I_{2n}]$,
respectively denote the collections of integer vectors $\mathbf{i}$
that lie in the interior,
the exterior, and the boundary parts of $R_n$, where $c_n$
is a suitably chosen sequence satisfying $c_n^{-1}+\lambda_n^{-1}c_n = \mathrm{o}(1)$.
It can be shown (cf. Section~\ref{sec64} below) that for $\beta>d+1/2 = 2.5$,
%
\begin{equation}\label{ex4-ord}
\sum_{k=1}^2 \sum
_{\mathbf{i}\in I_{kn}} \theta_n(\mathbf{i})^2 = \mathrm{o}(
\lambda_n)\quad  \mbox{and} \quad\sum_{\mathbf{i}\in I_{3n}}
\theta_n(\mathbf{i})^2 = \sigma^2_0
\lambda_n \bigl(1+\mathrm{o}(1)\bigr),
\end{equation}
where
\[
\sigma^2_0 = 16 B^2 \Biggl[B^2+
\sum_{k=1}^\infty \Biggl(\sum
_{j=k}^\infty b(j) \Biggr)^2 + \sum
_{k=1}^\infty \Biggl(\sum
_{j=k+1}^\infty b(j) \Biggr)^2 \Biggr].
\]
Hence, in this case, the contribution of the boundary part dominates
the other two terms, and the scaling is given by
${\lambda}_{n}\equiv\lambda_n^{d-1}$, which
is the $(d-1)$-dimensional Lebesgue measure of the boundary
of the sampling region $R_n$. Note that the rate of convergence
no longer depends on $\beta\in(d+1/2,\infty)$.
\end{exa}

The main reason why the edge effect dominates in the ND case as
highlighted by Example~\ref{ex42} can be explained by noting the form of
the constant $\sigma_0^2$ in (\ref{ex4-ord}). Although the
condition $A=0$
makes the sum of the $\alpha(\mathbf{i})$ over large open neighborhoods
of the origin small, sums of the $\alpha(\mathbf{i})$ over
half-planes, as determined by the $\theta_n(\mathbf{i})$
near the boundary of $R_n$ are not small.
As a result, the combined contribution
of these terms near the edge of $R_n$ determines the asymptotic
behavior of the sum $S_n$ for $\beta>d+1/2$.

The next result proves the CLT in presence of
non-trivial edge effects, for $\beta>d+1/2$.
The case $\beta=d+1/2$ will be
treated in Theorem~\ref{thm44} below.

\begin{theorem}\label{th43}
Let $\{Z(\cdot)\}$ be the linear process given by
(\ref{li-p}) such that conditions \textup{(C.1)}
hold
with
$\beta\in(d+1/2, \infty)$. Also suppose that $d\geq2$,
$A= 0$ and the following condition holds:
%
\begin{equation}\label{e-cond}
\lim_{\delta\downarrow0} \limsup_{n\to\infty}  \biggl| {
\lambda}_{n}^{-(d-1)} \sum_{\mathbf{i}\in[\partial R_n]^{\delta
{\lambda}_{n}}\cap\mathbb{Z}^d} \bigl|
\theta_n(\mathbf{i})\bigr|^2 - \sigma_{\mathrm{EE}}^2
\biggr| =0
\end{equation}
for some
$\sigma_{\mathrm{EE}}^2\in(0,\infty)$.
Then
%
\begin{equation}\label{ee-clt}
\lambda_n^{-(d-1)/2}[S_n - ES_n]
\mathop{\to}^d N\bigl(0, \sigma_{\mathrm{EE}}^2\bigr).
\end{equation}
\end{theorem}

Thus, it follows that under the conditions of the theorem, only
the $\theta(\mathbf{i})$ with indices $\mathbf{i}$ close to the
boundary of $R_n$ contribute to
the asymptotic variance of the sum. The
contribution of the $\theta(\mathbf{i})$ for the rest of $\mathbf{i}$-values
becomes asymptotically
negligible for $\beta>d+1/2$.
It can be shown that in Example~\ref{ex42}, the limiting variance $\sigma
_{\mathrm{EE}}^2$
is given by $\sigma^2_0$. Although we do not explicitly
state it, note that the boundary condition (\ref{bd-cd}) on $R_0$ is
implicit in the formulation
of (\ref{e-cond}). Also, note that this
edge-effect phenomenon in the ND case appears ONLY in
dimensions $d\geq2$.

Next, we consider the case where $\beta= d+1/2$. In this case, the
edge effect
may or may not have a non-trivial effect on the limit distribution,
depending on the growth rate of the slowly varying function $L(\cdot)$. More precisely,
we have the following results.

\begin{theorem}\label{thm44}
Let $\{Z(\cdot)\}$ be the linear process given by
(\ref{li-p}) such
that conditions \textup{(C.1)} holds with $\beta=d+1/2$ for some
$d\geq2$. Further suppose that $A=0$ and that (\ref{e-cond})
holds.
\begin{enumerate}[(ii)]
\item[(i)]
If $L(t) =\mathrm{o}(1)$ as $t \to\infty$, then (\ref{ee-clt})
holds.
\item[(ii)]  Suppose that condition \textup{(C.2)} holds and
${G^\dagger_{\infty}}\in L^2(\mathbb{R}^d)$ where ${G^\dagger
_{\infty}}$ is as in (\ref{cgi}).
\begin{enumerate}[(a)]
\item[(a)] If $L(t) =c_0(1+\mathrm{o}(1))$ as $t \to\infty$ for some $c_0\in
(0,\infty)$, then
\[
\lambda_n^{-(d-1)/2}[S_n - ES_n]
\mathop{\to}^d N \biggl(0, \sigma_{\mathrm{EE}}^2+
c_0^2 \int_{{\mathbb{R}}^d}
\bigl[{G^\dagger_{\infty}}(\mathbf{x})\bigr]^2\, \mathrm{d}
\mathbf{x} \biggr).
\]
\item[(b)]
If $L(t)^{-1} =\mathrm{o}(1)$ as $t \to\infty$,
then (\ref{clt-nd}) holds.
\end{enumerate}
\end{enumerate}
\end{theorem}

Theorem~\ref{thm44} shows that the edge effect is non-trivial under ND whenever
${\lambda}_{n}^{{3d-2\beta}} L^2({\lambda}_{n}) =\mathrm{O}(\lambda
_n^{(d-1)})$, that is, whenever the
contribution to $\sigma_n^2$ from the $\theta_n(\mathbf{i})^2$ near
the boundary is at
least as large as that from the remaining $\theta_n(\mathbf{i})^2$.
When the slowly varying function $L(\cdot)$ is bounded,
both ${\lambda}_{n}^{{3d-2\beta}} L^2({\lambda}_{n})$ and $\lambda
_n^{(d-1)}$ are of the same
order and the asymptotic variance depends on both $\sigma^2_{\mathrm
{EE}}$ and
the function ${G^\dagger_{\infty}}(\cdot)$ of (\ref{cgi}). On the
other\vspace*{1pt} hand, when the factors
${\lambda}_{n}^{{3d-2\beta}} L^2({\lambda}_{n})$ and $\lambda_n^{(d-1)}$
are not asymptotically equivalent, the scaling sequence and the
asymptotic variance of the centered sum are determined by the
dominant factor.

\begin{rem}\label{re45}
Note that in dimensions $d\geq2$, the
slowest possible growth rate of the variance of
the sum in the ND case is ${\lambda}_{n}^{(d-1)}$. This
may be contrasted with the
one dimensional ND case where the variance of the sum grows
at rate $\lambda_n^{[3-2\beta]}L({\lambda}_{n})^2$
which, in turn, can grow very slowly for $\beta$ close to $d+1/2=3/2$.
The main reason for this unusual behavior of the sum
in higher dimensions is the presence of the edge effect which
is not rate adaptive,
that is, it does not become asymptotically smaller even when the coefficients
$\alpha(\mathbf{i})$ or the function $\gamma(t)$ have a faster rate
of decay.
\end{rem}

\begin{rem}\label{re46}
For the one dimensional ND case, the variance
of the sum $\sigma^2_n$ does not necessarily go to infinity
for $\beta\geq3/2$ and hence, rate adaptivity of the
variance for $d=1$ is meaningful only when $\beta<3/2$, which is
covered by Theorem~\ref{th41}. It can be shown that
when $\beta=d+1/2$ and $d=1$, part (ii)(b)
of Theorem~\ref{thm44} holds. However, no analog of parts
(i) and (ii)(a) holds for $d=1$, as the edge effects
are asymptotically negligible in the one dimensional case. Also,
for $\beta\in(d+1/2,\infty)$, CLTs for the sum are not available
for $d=1$ in the ND case (as $\sigma_n^2\not\to\infty$), but they
are available
in dimensions $d\geq2$ (cf. Theorem~\ref{th43}).
\end{rem}

\section{Result under SRD}\label{sec5}
For completeness, we also give the result in the SRD case.
Suppose that $\beta\in[d,\infty)$ with $\int_0^\infty
t^{d-1}\gamma(t)\,\mathrm{d}t<\infty$. Then
it follows that $ A=\sum_{\mathbf{i}\in\mathbb{Z}^d} \alpha
(\mathbf{i})
\in\mathbb{R}$. If $A\neq0$, then the spatial process is SRD
and we have the following result.

\begin{theorem}\label{th51}
Let $\{Z(\cdot)\}$ be the linear process given by
(\ref{li-p}) such
that condition \textup{(C.1)} holds with
$\beta\in[d,\infty)$ and that $\int_0^\infty
t^{d-1}\gamma(t)\,\mathrm{d}t<\infty$ and $A\neq0$. Then, as $n\to\infty$,
%
\begin{equation}
\label{clt-srd} 
N_n^{-1/2}[S_n
-ES_n]  \mathop{\to}^d  N \bigl(0, A^2 \bigr).
\end{equation}
\end{theorem}

Thus, under the conditions of Theorem~\ref{th51}, the sum $S_n$ is
asymptotically normal
and the asymptotic variance grows at the standard rate, namely,
the square root of the sample size. Note that in this case,
we only assume condition (C.1), but not (C.2) or (C.3). The asymptotic
variance of the sum depends on the coefficients $\alpha(\cdot)$
only through the sum $A$, but not on the relative behavior of
the $\alpha(\cdot)$ and $\gamma(\cdot)$ at infinity.

\begin{rem}\label{rem52}
The asymptotic variance in the SRD case is determined
by the
$\alpha(\mathbf{i})$ for $\mathbf{i}$ in arbitrarily large compact
neighborhoods
of the origin, while in the PSD case, it is determined by the relative
behavior of $\alpha(\cdot)$ and $\gamma_t(\cdot)$ near infinity -- the
values of $\alpha(\mathbf{i})$ for any fixed compact neighborhood of
the origin
has no effect on the asymptotic variance.
In the ND case, the asymptotic variance depends on
the $\alpha(\mathbf{i})$ for both -- (i)
for smaller $\mathbf{i}$ through the condition $A=0$
and (ii) for large $\mathbf{i}$ through the relative
behavior of $g_t(\cdot)$ and $\gamma(\cdot)$ near infinity,
in absence of the edge-effect.
For $d\geq2$, in the ND case with non-trivial edge
effects, the asymptotic variance depends on the
$\alpha(\mathbf{i})$ for $\mathbf{i}$ in arbitrarily large
compact neighborhoods of the origin, as in the
SRD case, but not on the relative
behavior of $\alpha(\cdot)$ and $\gamma(\cdot)$.
\end{rem}

\section{Proofs}\label{sec6}

\subsection{Notation}\label{sec61}
Let $\mathcal{C}=[0,1)^d$ denote the unit cube in $\mathbb{R}^d$.
Let $L^2$ denote the collection of all square integrable functions
(with respect to the Lebesgue measure on $\mathbb{R}^d$) from
$\mathbb{R}^d$ to $\mathbb{R}$, and let $\|\cdot\|$ denote the $L^2$ norm,
that is, $\|f\|^2 = \int f^2(\mathbf{x})\, \mathrm{d}\mathbf{x}$, $f\in L^2$.
Let $\|\cdot\|_1$ and $\|\cdot\|_{\infty}$ , respectively
denote the $\ell^1$ and $\ell^{\infty}$-norms on $\mathbb{R}^d$,
that is,
for $(x_1,\ldots,x_d)'\in\mathbb{R}^d$,
\[
\bigl\|(x_1,\ldots,x_d)'\bigr\|_1 = \sum
_{i=1}^d |x_i|\quad \mbox{and}\quad
\bigl\|(x_1,\ldots,x_d)'\bigr\|_{\infty} = \max
\bigl\{ |x_i| \dvt  1\leq i\leq d\bigr\}.
\]
Recall that for any set $A\subset\mathbb{R}^d$ and $\delta\in
(0,\infty)$,
let $A^\delta=\{\mathbf{x}\in\mathbb{R}^d\dvt  \|\mathbf{x}-\mathbf
{y}\|\leq\delta$
for some $\mathbf{y}\in A\}$ denote the $\delta$-enlargement of $A$.
Similarly, define the set $A^{-\delta} = \{\mathbf{x}\in A\dvt  B(\mathbf
{x};\delta)\subset A\}$
where $ B(\mathbf{x};\delta) =\{\mathbf{y}\in\mathbb{R}^d\dvt  \|
\mathbf{x}-\mathbf{y}\|<\delta\}$.
Let $\partial A$ denote the boundary of $A$.
Let $C, C(\cdot)$ denote generic constants that do not depend on $n$.
Unless otherwise specified, all limits (including those in
the order symbols) are taken by letting
$n\to\infty$. Also, for notational simplicity, we set $\mu= 0$
for the rest of this section (except in cases where
there is a chance of confusion).

\subsection{Auxiliary results}\label{sec62}

Here we prove a general version of the CLT for spatial linear processes
without structural conditions on the coefficients
$\alpha(\mathbf{i})$ and the sampling regions. This result
forms the basis for proving the results from Sections~\ref{sec3}--\ref{sec5},
and may be of independent interest.

\begin{theorem}\label{thm61}
Let $\{Z(\mathbf{j})\dvt  \mathbf{j}\in\mathbb
{Z}^d\}$ be the spatial
linear process in (\ref{li-p}) with $\mu=0$.
Let $\Lambda_n$ be a finite subset of
$\mathbb{Z}^d$ such that $|\Lambda_n|\to\infty$ as $n\to\infty$.
Let $S_n=\sum_{\mathbf{j}\in
\Lambda_n} Z(\mathbf{j})$ and $\sigma_n^2 = \sum_{\mathbf{i}\in
\mathbb{Z}^d} \theta_n(\mathbf{i})^2$,
where $\theta_{n}(\mathbf{i}) = \sum_{\mathbf{j}\in\Lambda_n}
\alpha(\mathbf{j}-\mathbf{i})$,
$\mathbf{i}\in\mathbb{Z}^d$. Suppose that as $n\to\infty$,
%
\begin{equation}\label{clt-cond}
\frac{1}{\sigma_n} 
+ \frac{\max\{| \theta_{n}(\mathbf{i})| \dvt  \mathbf{i}\in\mathbb
{Z}^d\}}{\sigma_n} \to0.
\end{equation}
Then $S_n/\sigma_n \mathop{\to}^d N(0,1)$ as $n\to\infty$.
\end{theorem}

\begin{pf}
First, we shall show that there exists
a sequence of integers $m_n\to\infty$
such that
%
\begin{equation}\label{a-neg1}
\sigma_n^{-2}\sum_{\|\mathbf{i}\|> m_n}
\theta_n(\mathbf{i})^2 = \mathrm{o}(1).
\end{equation}
To that end, note that $|\theta_n(\mathbf{i})| \leq[ \sum_{\mathbf
{j}\in(\Lambda_n -\mathbf{i})}\alpha(\mathbf{j})^2]^{1/2} |\Lambda
_n|^{1/2}$
for all $\mathbf{i}\in\mathbb{Z}^d$. Since $\sum_{\mathbf{i}\in
\mathbb{Z}^d}\alpha(\mathbf{i})^2<\infty$,
there exists $m_{1n}\to\infty$ such that
$
\sigma_n^{-2}|\Lambda_n|^2 \sum_{\|\mathbf{i}\|>m_{1n}}\alpha
(\mathbf{i})^2 = \mathrm{o}(1)
$.
Define $m_n = \max\{\|\mathbf{j}\|\dvt  \mathbf{j}\in\Lambda_n\} +
m_{1n}$. Then,
it follows that
\[
\sigma_n^{-2}\sum_{\|\mathbf{i}\|> m_n}
\theta_n(\mathbf{i})^2 \leq \sigma_n^{-2}
|\Lambda_n| \sum_{\|\mathbf{i}\|> m_n} \sum
_{\mathbf
{j}\in\Lambda_n} \alpha(\mathbf{j}-\mathbf{i})^2 \leq
\sigma_n^{-2} |\Lambda_n|^2 \sum
_{\|\mathbf{i}\|>m_{1n}} \alpha (\mathbf{i})^2 = \mathrm{o}(1),
\]
proving (\ref{a-neg1}).

Next, define $U_n = \{\mathbf{i}\in\mathbb{Z}^d\dvt  \|\mathbf{i}\|\leq
m_n\}$
and $\bar{U}_n = \{\mathbf{i}\in\mathbb{Z}^d\dvt  \|\mathbf{i}\|> m_n\}$. Define
$\tilde{S}_n = \sum_{\mathbf{i}\in U_n} \theta_n(\mathbf
{i})\varepsilon(\mathbf{i})$ and
$\tilde{\sigma}_n^2 = \sum_{\mathbf{i}\in U_n} \theta_n(\mathbf{i})^2$.
Then, by (\ref{a-neg1}), $\sigma_n^{-2}[\sigma_n^2 - \tilde{\sigma
}_n^2] = \mathrm{o}(1)$
and hence,
%
\begin{eqnarray}
\sigma_n^{-1}S_n - \tilde{
\sigma}_n^{-1} \tilde{S}_n &=&
\sigma_n^{-1} \biggl[\sum_{\mathbf{i}\in\bar{U}_n}
\theta _n(\mathbf{i}) \varepsilon(\mathbf{i}) \biggr] +
\frac{[\tilde{\sigma}_n - \sigma_n]}{\sigma_n \tilde{\sigma
}_n} \tilde{S}_n
\nonumber
\\[-8pt]
\label{a-neg2}
\\[-8pt]
\nonumber
&=&
\mathrm{o}_p(1), 
\end{eqnarray}
provided $\tilde{\sigma}_n^{-1}\tilde{S}_n= \mathrm{O}_p(1)$.
Hence, it is enough to show that $\tilde{\sigma}_n^{-1}\tilde
{S}_n\mathop{\to}^d N(0,1)$
as $n\to\infty$.

By Lindeberg's CLT, this would follow if
for all $\delta>0$,
%
\begin{equation}
\label{LC}
\lim_{n\to\infty}\sum_{\mathbf{i}\in U_n}
E Y_n(\mathbf {i})^2\mathbh{1}\bigl(\bigl|Y_n(
\mathbf{i})\bigr|> \delta\bigr) = 0,
\end{equation}
where $Y_n(\mathbf{i}) = \varepsilon(\mathbf{i}) \theta_n(\mathbf
{i})/\tilde{\sigma}_n$, $\mathbf{i}\in U_n$.
Now, by uniform integrability of $\{\varepsilon(\mathbf{i})^2\dvt  \mathbf
{i}\in\mathbb{Z}^d\}$
and (\ref{clt-cond}), for any $\delta>0$,
\begin{eqnarray*}
\sum_{\mathbf{i}\in U_n} E Y_n(
\mathbf{i})^2\mathbh {1}\bigl(\bigl|Y_n(\mathbf{i})\bigr|> \delta
\bigr)
&=&  \tilde{\sigma}_n^{-2} \sum
_{\mathbf{i}\in U_n} \theta _n(\mathbf{i})^2 E
\varepsilon(\mathbf{i})^2\mathbh{1}\bigl(\bigl|\theta_n(
\mathbf{i}) \varepsilon (\mathbf{i})\bigr|> \delta\tilde{\sigma}_n\bigr)
\\
&\leq& \max_{\mathbf{i}\in\mathbb{Z}^d} E \varepsilon(\mathbf{i})^2
\mathbh{1} \biggl(\bigl|\varepsilon(\mathbf{i})\bigr|> \delta\frac{\tilde{\sigma}_n}{\max_{\mathbf{j}\in\mathbb{Z}^d}
|\theta_n(\mathbf{j})|} \biggr)
\\
&= &\mathrm{o}(1).
\end{eqnarray*}
Hence, (\ref{LC}) holds and the result is proved.
\end{pf}

\begin{cor}\label{co62}
Let $\Lambda_n$ be as in Theorem~\ref{thm61}. Then,
(\ref{clt-cond}) holds if either of the
following two conditions holds:
\begin{longlist}[(ii)]
\item[(i)]
$\max\{| \theta_{n}(\mathbf{i})| \dvt  \mathbf{i}\in\mathbb{Z}^d\}=
\mathrm{O}(1)$ and
${\sigma_n}^2 \to\infty$ as $n\to\infty$.
\item[(ii)] $\liminf_{n\to\infty} \sigma_n^2/|\Lambda_n| >0$.
\end{longlist}
\end{cor}

\begin{pf}
Sufficiency of (i) for (\ref{clt-cond}) is trivial.
Consider (ii).
Fix a sequence $\{p_n\}\subset(0,\infty)$ such that
$p_n^{-1}+ |\Lambda_n|^{-1/d} p_n=\mathrm{o}(1)$.
Then, by the Cauchy--Schwarz inequality (and (ii)),
\begin{eqnarray*}
&&\max_{\mathbf{j}\in\mathbb{Z}^d}\frac{| \theta_{n}(\mathbf
{j})|}{\sigma_n}\\
&&\quad\leq\max_{\mathbf{j}\in\mathbb{Z}^d}
\sigma_n^{-1} \biggl[\sum_{\|\mathbf{i}- \mathbf{j}\| \leq p_n, \mathbf{i}\in\Lambda_n}
\bigl|\alpha(\mathbf{i}-\mathbf{j})\bigr| + \sum_{\|\mathbf{i}- \mathbf{j}\| > p_n, \mathbf{i}\in\Lambda_n} \bigl|
\alpha(\mathbf{i}-\mathbf{j})\bigr|  \biggr]
\\
&&\quad\leq \sigma_n^{-1} \biggl[ C(d) p_n^{d/2}
\biggl(\sum_{\mathbf{i}\in\mathbb{Z}^d} \alpha(\mathbf {i})^2
\biggr)^{1/2} + |\Lambda_n|^{1/2} \biggl(\sum
_{\|\mathbf{i}\|> p_n} \alpha(\mathbf {i})^2
\biggr)^{1/2} \biggr]
\\
&&\quad= \frac{|\Lambda_n|^{1/2}}{\sigma_n} \bigl[\mathrm{O} \bigl(|\Lambda _n|^{-1/2}
p_n^{d/2} \bigr) + \mathrm{o}(1) \bigr] = \mathrm{o}(1).
\end{eqnarray*}
This completes the proof of the corollary.
\end{pf}

\subsection{Proofs of the results from Section~\texorpdfstring{\protect\ref{sec3}}{3}}\label{sec63}

\begin{pf*}{Proof of Proposition~\ref{pr3.1}}
Note that by definition,
$|g_t(\mathbf{z})| \leq1$ for all $\|\mathbf{z}\|=1$ and $t>0$. For any
$t>0$ and for any
$\mathbf{x}\neq\mathbf{0}$, writing $\mathbf{x}=r\mathbf{z}$ with
$\|\mathbf{z}\|=1$
and $r=\|\mathbf{x}\|$, we have
%
\begin{equation}\label{gt-bd}
\bigl|g_t(\mathbf{x})\bigr| = \bigl|g_{tr}(\mathbf{z})\bigr|
\frac{\gamma(tr)}{\gamma
(t)} \leq\frac{\gamma(tr)}{\gamma(t)} = \|x\|^{-\beta}
\frac{L(t\|\mathbf{x}\|)}{L(t)}.
\end{equation}
Next, using condition (C.2) and a subsequence argument
(cf. page 92, Athreya and Lahiri \cite{AthLah06}), we have
%
\begin{equation}\label{gi-bd}
\bigl|{g}_{\infty}(\mathbf{x})\bigr| \leq \|x\|^{-\beta} \qquad \mbox{almost
everywhere ($m$)}, \mathbf{x}\neq\mathbf{0},
\end{equation}
where $m$ is the Lebesgue measure on $\mathbb{R}^d$.
Hence,
$
\int_{R_0} \|{g}_{\infty}(\mathbf{y}-\mathbf{x})\|\, \mathrm{d}\mathbf{y}\leq
\int_{R_0} \| \mathbf{y}-\mathbf{x}\|^{-\beta}\mathbh{1}(\mathbf
{y}\neq\mathbf{x})\, \mathrm{d}\mathbf{y}<\infty
$
for all $\mathbf{x}\in\mathbb{R}^d$ and for all $\beta\in(0,d)$.
This completes the proof.
\end{pf*}

\begin{pf*}{Proof of Theorem~\ref{th32}}
First, we shall show that ${G}_{\infty}\in
L^2(\mathbb{R}^d)$. Note that
$R_0\subset B(\mathbf{0},2^{-1}\sqrt{d})$ and hence,
by (\ref{gi-bd}),
$|{G}_{\infty}(\mathbf{x})|\leq\int_{\|\mathbf{y}\|\leq2\sqrt
{d}} |{g}_{\infty}(\mathbf{y})|\,\mathrm{d}\mathbf{y}\leq C(d,\beta) $ for all
$\|\mathbf{x}\|^2\leq d$
while $|{G}_{\infty}(\mathbf{x})|\leq\int_{\|\mathbf{y}-\mathbf
{x}\|\leq\sqrt{d}/2} |{g}_{\infty}(\mathbf{y})|\,\mathrm{d}\mathbf{x}\leq
C(d,\beta) \|\mathbf{x}\|^{-\beta}$ for all
$\|\mathbf{x}\|^2> d$, implying that ${G}_{\infty}\in L^2(\mathbb{R}^d)$.

Next, we apply Corollary~\ref{co62} to establish Theorem~\ref{th32}. Note that by
condition (C.2), there exists a sequence
${\eta}_{n}\downarrow0$ such that
%
\begin{equation}\label{c22}
\int_{\{2\|\mathbf{x}\|>{\eta}_{n}\}} \bigl|{g}_{\infty}(\mathbf {x})
-g_{{\lambda}_{n}}(\mathbf{x}) \bigr|^2 \,\mathrm{d}\mathbf{x} = \mathrm{o}(1)\qquad \mbox{as } n\to
\infty.
\end{equation}
W.l.o.g, suppose that ${\lambda}_{n}{\eta}_{n}\gg{\lambda
}_{n}^\delta$ for some $\delta\in(0,1/2)$. Next, for $\mathbf{i}\in
\mathbb{Z}^d$, write
\begin{eqnarray*}
\theta_n(\mathbf{i}) &=& \sum_{\mathbf{j}\in(R_n-\mathbf{i})\cap
\mathbb{Z}^d, \|\mathbf{j}\| > {\lambda}_{n}{\eta}_{n}}
\alpha (\mathbf{j}) + \sum_{\mathbf{j}\in(R_n-\mathbf{i})\cap\mathbb{Z}^d, \|\mathbf
{j}\| \leq{\lambda}_{n}{\eta}_{n}} \alpha(\mathbf{j})
\\
&\equiv& \theta_{1n}(\mathbf{i}) +\theta_{2n}(\mathbf{i}),\qquad
\mbox{say.}
\end{eqnarray*}

We shall first show that the contribution
from the $\theta_{2n}(\mathbf{i})$-terms to $\sigma_n^2$ is
negligible. To that end, note that by definition,
$\theta_{2n}(\mathbf{i})=0$ for all $\mathbf{i}\notin[-2{\lambda
}_{n}, 2{\lambda}_{n}]^d$. Hence,
\begin{eqnarray}
\sum_{\mathbf{i}\in\mathbb{Z}^d} \theta_{2n}(
\mathbf{i})^2 & \leq& \sum_{\mathbf{i}\in [-2{\lambda}_{n}, 2{\lambda}_{n}]^d} \biggl(
\sum_{\|\mathbf{j}\| \leq{\lambda}_{n}{\eta}_{n}} \bigl|\alpha(\mathbf {j})\bigr|
\biggr)^2
\nonumber\\
&\leq& (4{\lambda}_{n})^d \biggl( \sum
_{\|\mathbf{j}\| \leq{\lambda
}_{n}{\eta}_{n}} \bigl|\alpha(\mathbf{j})\bigr| \biggr)^2
\nonumber
\\[-8pt]
\label{t2n}\\[-8pt]
\nonumber
&\leq& C(d) {\lambda}_{n}^d \bigl( [{
\lambda}_{n}{\eta }_{n}]^{d-\beta} L({
\lambda}_{n}{\eta}_{n}) \bigr)^2
\\
\nonumber
&=& \mathrm{o} \bigl({\lambda}_{n}^{3d-2\beta}L^2({
\lambda}_{n}) \bigr).
\end{eqnarray}

Next, consider the $\theta_{1n}(\mathbf{i})$-terms. Note that by
definition, for any $t>0$ and $\mathbf{k}\in\mathbb{Z}^d$,
$g_{t}(\mathbf{x}) = g_t({\mathbf{k}})$
for all $\mathbf{x}\in{t}^{-1}(\mathbf{k}+\mathcal{C})$. Hence,
\begin{eqnarray*}
\theta_{1n}(\mathbf{i}) 
&=& \sum
_{\mathbf{j}\in (R_n - \mathbf{i})\cap\mathbb{Z}^d, \|
\mathbf{j}\|>{\lambda}_{n}{\eta}_{n}} \alpha(\mathbf{j})
\\
&=& \sum_{\mathbf{j}\in (R_n - \mathbf{i})\cap\mathbb{Z}^d, \|
\mathbf{j}\|>{\lambda}_{n}{\eta}_{n}} g_{{\lambda}_{n}}(\mathbf {j}/{
\lambda}_{n})\gamma({\lambda}_{n})
\\
&=& \gamma({\lambda}_{n}) {\lambda}_{n}^d \int
g_{{\lambda
}_{n}}(\mathbf{x}) \mathbh{1} \bigl(\mathbf{x}\in \bigl\llbracket {R_0 - {\lambda}_{n}^{-1}{\mathbf{i}}}
\bigr\rrbracket_n \bigr) \,\mathrm{d}\mathbf{x} 
\\
&\equiv& \gamma({\lambda}_{n}) {\lambda}_{n}^d
G_n(\mathbf {i}/{\lambda}_{n}),\qquad\mbox{say,}
\end{eqnarray*}
where
$\llbracket {R_0 -\mathbf{y}}\rrbracket_n =
\bigcup\{ {\lambda}_{n}^{-1}(\mathbf{k}+\mathcal{C}) \dvt
\mathbf{k}\in\mathbb{Z}^d,
\|\mathbf{k}\|>{\lambda}_{n}{\eta}_{n},
\frac{\mathbf{k}}{{\lambda}_{n}}\in R_0 -\mathbf{y}\}$
and $G_n(\mathbf{y}) = \int g_{{\lambda}_{n}} (\mathbf{x})\times \mathbh
{1}(\mathbf{x}\in
\llbracket{R_0 - \mathbf{y}}\rrbracket_n)\, \mathrm{d}\mathbf{x}$, $\mathbf{y}\in
\mathbb{R}^d$.
Note that
\begin{eqnarray*}
\sigma_{1n}^2&\equiv& \sum_{\mathbf{i}\in\mathbb{Z}^d}
\theta _{1n}(\mathbf{i})^2
\\
&=& \gamma({\lambda}_{n})^2 {\lambda}_{n}^{2d}
\sum_{\mathbf{i}\in
\mathbb{Z}^d} G_n^2(
\mathbf{i}/{\lambda}_{n})
\\
&=& \gamma({\lambda}_{n})^2 {\lambda}_{n}^{3d}
\int\tilde {G}_n^2(\mathbf{x})\,\mathrm{d}\mathbf{x},\qquad\mbox{say,}
\end{eqnarray*}
where $\tilde{G}_n(\mathbf{x}) = \sum_{\mathbf{i}\in\mathbb{Z}^d}
G_n(\mathbf{i}/{\lambda}_{n})
\mathbh{1} (\mathbf{x}\in{\lambda}_{n}^{-1}(\mathbf
{i}+\mathcal{C}) )$, $\mathbf{x}\in\mathbb{R}^d$.
We shall now show that
%
\begin{equation}\label{t1-cg}
\|\tilde{G}_n - {G}_{\infty}\|_2\to0\qquad \mbox{as } n
\to\infty.
\end{equation}

To that end, 
write $a_n=\sqrt{d}\lambda_n^{-1}$ and
note that
$\sup\{ \|{\lambda}_{n}^{-1}{\mathbf{i}} -\mathbf{z}\| \dvt  \mathbf
{z}\in{\lambda}_{n}^{-1}(\mathbf{i}+\mathcal{C})\}
\leq a_n$ and that
$
\llbracket{R_0 - {\lambda}_{n}^{-1}{\mathbf{i}}}\rrbracket_n \subset[R_0^{a_n}
- {\lambda}_{n}^{-1}{\mathbf{i}}]\cap\{\mathbf{x}\in\mathbb{R}^d
\dvt  \|\mathbf{x}\|> {\eta}_{n}- a_n\}
$. Write $g_t^{(n)}(\mathbf{x}) = g_t(\mathbf{x})\mathbh{1}( \|
\mathbf{x}\|> {\eta}_{n}- a_n)$,
$t\in(0,\infty]$ and $\mathbf{x}\in\mathbb{R}^d$. Then,
%
\begin{eqnarray}
&& \sum_{\mathbf{i}\in\mathbb{Z}^d} \int
_{{\lambda
}_{n}^{-1}(\mathbf{i}+\mathcal{C})} \biggl( \int \bigl[ g_{{\lambda}_{n}}(\mathbf{x}) -
{g}_{\infty}(\mathbf {x})\bigr]\mathbh{1} \bigl(\mathbf{x}\in \bigl\llbracket  {R_0 - {\lambda}_{n}^{-1}{\mathbf{i}}}
\bigr\rrbracket_n  \bigr)\,\mathrm{d}\mathbf {x} \biggr)^2\, \mathrm{d}\mathbf{z}\nonumber
\\
&&\quad \leq \sum_{\mathbf{i}\in\mathbb{Z}^d} \int_{{\lambda}_{n}^{-1}(\mathbf
{i}+\mathcal{C})}
\biggl( \int \bigl| g_{{\lambda}_{n}}(\mathbf{x}) - {g}_{\infty}(\mathbf {x}) \bigr|
\mathbh{1} \bigl(\mathbf{x}\in R_0^{2a_n} - \mathbf{z} \bigr)
\mathbh{1} \bigl(\|\mathbf{x}\|> {\eta}_{n}- a_n \bigr) \,\mathrm{d}\mathbf{x}
\biggr)^2 \,\mathrm{d}\mathbf{z}\nonumber
\\
&& \quad=
\int\!\!\!\int\!\!\!\int \bigl| g_{{\lambda}_{n}}^{(n)}(\mathbf{x}-\mathbf{z}) -
{g}_{\infty}^{(n)}(\mathbf{x}-\mathbf{z}) \bigr| \bigl| g_{{\lambda}_{n}}^{(n)}
(\mathbf{y}-\mathbf{z}) - {g}_{\infty
}^{(n)}(\mathbf{y}-
\mathbf{z}) \bigr|\nonumber
\\
&&\hspace*{13pt}\quad\qquad{}\times\mathbh{1} \bigl(\mathbf{x}\in R_0^{2a_n}
\bigr) \mathbh{1} \bigl(\mathbf{y}\in R_0^{2a_n} \bigr)\, \mathrm{d}
\mathbf{x}\, \mathrm{d}\mathbf{y} \,\mathrm{d}\mathbf{z}\nonumber
\\
\label{diff1}&&\quad\leq \int\!\!\!\int \biggl(\int \bigl| g_{{\lambda}_{n}}^{(n)}(\mathbf{x}-
\mathbf{z}) - {g}_{\infty}^{(n)}(\mathbf{x}-\mathbf{z})
\bigr|^2 \,\mathrm{d}\mathbf{z} \biggr)^{1/2} \\
&&\hspace*{6pt}\quad\qquad{}\times\biggl(\int \bigl|
g_{{\lambda}_{n}}^{(n)}(\mathbf{y}-\mathbf{z}) - {g}_{\infty
}^{(n)}(
\mathbf{y}-\mathbf{z}) \bigr|^2 \,\mathrm{d}\mathbf{z} \biggr)^{1/2}
\nonumber\\
&&\hspace*{6pt}\quad\qquad{}\times \mathbh{1} \bigl(\mathbf{x}\in R_0^{2a_n}
\bigr) \mathbh{1} \bigl(\mathbf{y}\in R_0^{2a_n} \bigr)\, \mathrm{d}
\mathbf{x}\, \mathrm{d}\mathbf{y}\nonumber
\\
&& \quad= \bigl\|g_{{\lambda}_{n}}^{(n)} - {g}_{\infty}^{(n)}
\bigr\|_2^2 \biggl[\int \mathbh{1} \bigl(\mathbf{x}\in
R_0^{2a_n} \bigr)\,\mathrm{d}\mathbf{x} \biggr]^2\nonumber
\\
&&\quad=\mathrm{o}(1).\nonumber
\end{eqnarray}
Next, note that
the symmetric difference of the sets
$ \llbracket R_0 - {\lambda}_{n}^{-1}{\mathbf{i}} \rrbracket_n$ and $R_0 - \mathbf{z}$
is contained in $(\partial R_0 -\mathbf{z})^{2a_n}
= (\partial R_0)^{2a_n} -\mathbf{z}$
for all $\mathbf{z}\in{\lambda}_{n}^{-1}(\mathbf{i}+\mathcal{C})$
and for all $\mathbf{i}\in\mathbb{Z}^d$
with $\|\mathbf{i}\|>C(d)\lambda_n$, while it is contained in
$[(\partial R_0)^{2a_n} -\mathbf{z}]\cup\{\mathbf{x}\dvt \|\mathbf{x}\|
\leq{\eta}_{n}+a_n\}$
for all
$\mathbf{z}\in{\lambda}_{n}^{-1}(\mathbf{i}+\mathcal{C})$, for $\|
\mathbf{i}\|\leq C(d)\lambda_n$.
Now, using the set inclusion relations given above (for the first inequality),
arguments similar to~(\ref{diff1}) above (for the first term of the
last inequality),
and
the bounded convergence
theorem and the regularity conditions on the boundary of $R_0$ (that
the $d$-dimensional Lebesgue measure of $\partial R_0$ is zero)
(for the second and the third terms in the last inequality),
we have
%
\begin{eqnarray}
&& \sum_{\mathbf{i}\in\mathbb{Z}^d} \int_{{\lambda
}_{n}^{-1}(\mathbf{i}+\mathcal{C})} \biggl(
\int {g}_{\infty}(\mathbf{x}) \bigl[\mathbh{1} \bigl(\mathbf{x}\in
 \bigl\llbracket  {R_0 - {\lambda}_{n}^{-1}{\mathbf{i}}}
\bigr\rrbracket_n  \bigr) - \mathbh{1} (\mathbf{x}\in R_0 - \mathbf{z} )
\bigr] \,\mathrm{d}\mathbf{x} \biggr)^2 \,\mathrm{d}\mathbf{z}
\nonumber\\
&& \quad \leq \sum_{\|\mathbf{i}\|>C(d)\lambda_n} \int_{{\lambda
}_{n}^{-1}(\mathbf{i}+\mathcal{C})}
\biggl\{ \int\bigl| {g}_{\infty}(\mathbf{x})\bigr| \mathbh{1} \bigl(\mathbf{x}\in (
\partial R_0)^{2a_n} -\mathbf{z} \bigr)\,\mathrm{d}\mathbf{x} \biggr
\}^2 \,\mathrm{d}\mathbf{z}
\nonumber\\
&& \qquad{}+ \sum_{\|\mathbf{i}\|\leq C(d)\lambda_n} \int_{{\lambda
}_{n}^{-1}(\mathbf{i}+\mathcal{C})}
\biggl[ \int\bigl| {g}_{\infty}(\mathbf{x})\bigr| \bigl\{ \mathbh{1} \bigl(\mathbf{x}
\in (\partial R_0)^{2a_n} -\mathbf{z} \bigr) + \mathbh{1} \bigl(\|
\mathbf{x}\| \leq{\eta}_{n}+ a_n \bigr) \bigr\} \,\mathrm{d}\mathbf{x}
\biggr]^2\, \mathrm{d}\mathbf{z}
\nonumber\\
\label{diff2}
&& \quad \leq \bigl\|{g}_{\infty}(\mathbf{x})\mathbh{1} \bigl(\|\mathbf{x}\| > C(d)
\bigr) \bigr\|_2^2 \biggl\{\int\mathbh{1} \bigl(\mathbf{x}\in (
\partial R_0)^{2a_n} \bigr)\,\mathrm{d}\mathbf{x} \biggr\}^2
\\
&&\qquad {}+ 2 \int_{\|\mathbf{z}\|\leq C(d) } \biggl\{\int \bigl|
{g}_{\infty}(\mathbf{x})\bigr| \mathbh{1} \bigl( \mathbf{x}\in\partial
R_0^{2a_n} - \mathbf{z} \bigr)\,\mathrm{d}\mathbf{x} \biggr
\}^2\,\mathrm{d}\mathbf{z}\nonumber
\\
&& \qquad {}+ C(d) \biggl\{\int \bigl| {g}_{\infty}(\mathbf{x})\bigr|
\mathbh{1} \bigl(\|\mathbf {x}\| \leq{\eta}_{n}+a_n \bigr)\,\mathrm{d}
\mathbf{x} \biggr\}^2\nonumber
\\
&& \quad=\mathrm{o}(1),\nonumber
\end{eqnarray}
as $\int|{g}_{\infty}(\mathbf{x})|\mathbh{1}(\|\mathbf{x}\|\leq C)
\,\mathrm{d}\mathbf{x}+
\int{g}_{\infty}(\mathbf{x})^2\mathbh{1}(\|\mathbf{x}\|\geq C)
\,\mathrm{d}\mathbf{x}
<\infty$ for all $C\in(0,\infty)$ and for all $\beta\in(d/2,d)$.
This completes the proof of (\ref{t1-cg}).
Note that (\ref{t1-cg}) implies that
$\int\tilde{G}_n^2(\mathbf{x})\,\mathrm{d}\mathbf{x}\to\int{G}_{\infty
}^2(\mathbf{x})\,\mathrm{d}\mathbf{x}$
as $n\to\infty$. Hence, by Corollary~\ref{co62}(ii), Theorem~\ref{th32} follows.
\end{pf*}

\subsection{Proofs of the results from Section~\texorpdfstring{\protect\ref{sec4}}{4}}\label{sec64}

\begin{pf*}{Proof of Theorem~\ref{th41}}
By Corollary~\ref{co62}(i),
it is enough to show that ${G^\dagger_{\infty}}\in L^2$
and
%
\begin{equation}\label{si-nd}
\sigma_n^2 = \biggl[ {\lambda}_{n}^{3d-2\beta}L({
\lambda}_{n})^2 \int{G^\dagger_{\infty}}(
\mathbf{x})^2 \,\mathrm{d}\mathbf{x} \biggr]\bigl(1+\mathrm{o}(1)\bigr).
\end{equation}
Note that by (\ref{gi-bd}),
\[
\bigl|{G^\dagger_{\infty}}(\mathbf{x})\bigr| \leq C(d,\beta) \cases{
d(\mathbf{x}, R_0)^{d-\beta}, &\quad$\mbox{if }\mathbf{x}\in[R_0\cup \partial R_0]^c
\mbox{ and } \|\mathbf{x}\| \leq1$,
\vspace*{2pt}\cr
\|\mathbf{x}\|^{-\beta}, & \quad$\mbox{if }\|\mathbf{x}\| >1$}
\]
and $|{G^\dagger_{\infty}}(\mathbf{x})| \leq C(d,\beta) d(\mathbf
{x}, R_0^c)^{d-\beta}$ for $\mathbf{x}\in R_0$.
Since $\beta\in(d, d+1/2)$, by the boundary condition (C.3), it
follows that
${G^\dagger_{\infty}}\in L^2$.

Next, consider (\ref{si-nd}). For $\beta\in(d, d+1/2)$, by condition
(C.2), (\ref{gt-bd}) and
(\ref{gi-bd}), there exists a sequence
${\eta}_{n}\downarrow0$ such that
%
\begin{equation}\label{c2-et}
\int_{\{2\|\mathbf{x}\|>{\eta}_{n}\}} \bigl|{g}_{\infty}(\mathbf{x})
-g_{{\lambda}_{n}}(\mathbf{x}) \bigr|^p \,\mathrm{d}\mathbf{x} = \mathrm{o}(1)\qquad\mbox{as } n\to
\infty,
\end{equation}
for $p=1,2$. For $p=1$, this follows directly from (C.2), as $\beta>d$.
As for $p=2$, we use the trivial bound ``$f(x)^2 \leq|f(x)|
\sup_{x\in B}\{ |f(x)|\}$
for a function $f\dvtx B\to\mathbb{R}$'' in conjunction with the $p=1$ relation
and the bounds (\ref{gt-bd}) and (\ref{gi-bd}), which may require
replacing the ${\eta}_{n}$ for
the $p=1$ case by a possibly coarser sequence that still decreases to zero.
Hence, (\ref{c2-et}) holds for both $p=1,2$.

W.l.o.g., suppose that ${\eta}_{n}\gg{\lambda}_{n}^{-\delta}$ (i.e.,
${\lambda}_{n}^\delta{\eta}_{n}\to\infty$)
for some $\delta\in(0,1)$.
Let $t_n = {\lambda}_{n}{\eta}_{n}$ and let $u_n = {\lambda
}_{n}^{[2d-2\beta+1]/2}$. Then $u_n^{-1}+t_n^{-1}
= \mathrm{o}(1)$ and
$
u_n = \mathrm{o} ({\lambda}_{n}^{[2d-2\beta+1]}L_n({\lambda}_{n})^2  )
$.

Also, define
%
\begin{eqnarray}
V_{1n} &=& \bigl\{\mathbf{i}\in\mathbb{Z}^d\dvt  B(
\mathbf{i};t_n)\subset R_n\bigr\},\nonumber
\\
\label{v-s} V_{2n} &=& \bigl\{\mathbf{i}\in\mathbb{Z}^d\dvt  B(
\mathbf{i}; t_n)\subset R_n^c\bigr\},
\\
V_{3n} & =& \mathbb{Z}^d\setminus[V_{1n}\cup
V_{2n}].\nonumber
\end{eqnarray}

First, consider the sum of $\theta_n(\mathbf{i})^2$ for $\mathbf
{i}\in V_{3n}$. Note that
$V_{3n}\subset\{\mathbf{i}\in\mathbb{Z}^d\dvt  \mathbf{i}\in(\partial
R_n)^{t_n}\}$.
By (C.3) (with $f\equiv1$), $\nu((\partial R_0)^\varepsilon)
= \mathrm{O}(\varepsilon)$ as $\varepsilon\to0$.
This implies $\nu( (\partial R_n)^{u_n}) = \mathrm{O}({\lambda}_{n}^{d-1}u_n)$
and\vspace*{-2pt} hence,
%
\begin{equation}\label{i-v3n}
\sum_{\mathbf{i}\in\partial R_n^{u_n}\cap\mathbb{Z}^d} \bigl|\theta _n(
\mathbf{i})\bigr|^2 \leq  C(d) {\lambda}_{n}^{d-1}u_n
\biggl[\sum_{\mathbf{i}\in\mathbb{Z}^d} \bigl|\alpha(\mathbf{i})\bigr|
\biggr]^2 
= \mathrm{o} \bigl({
\lambda}_{n}^{[3d -2\beta]}L({\lambda}_{n})^2
\bigr). 
\end{equation}
If $u_n >t_n$, then this shows that $\sum_{\mathbf{i}\in V_{3n}}
\theta_n(\mathbf{i})^2
= \mathrm{o} ({\lambda}_{n}^{[3d -2\beta]}L({\lambda}_{n})^2)$.
Hence, w.l.o.g., suppose that $u_n\leq t_n$.
Let $v_n(\mathbf{i})= d(\mathbf{i}, \partial R_n)$, $\mathbf{i}\in
[(\partial R_n)^{u_n}]^c$.
By the condition $A=0$, uniformly in $\mathbf{i}\notin
(\partial R_n)^{u_n}$, we\vspace*{-2pt} have
\[
\bigl|\theta_n(\mathbf{i})\bigr|\leq\sum_{\|\mathbf{l}\|>v_n(\mathbf{i})} \bigl|
\alpha(\mathbf{l})\bigr| \leq \sum_{\|\mathbf{l}\|>v_n(\mathbf{i}) } \|\mathbf{l}
\|^{-\beta
} L\bigl(\|\mathbf{l}\|\bigr) \leq C(d,\beta) v_n(
\mathbf{i})^{d-\beta}L\bigl(v_n(\mathbf{i})\bigr).
\]
Hence, by condition (C.3), it follows\vspace*{-2pt} that
%
\begin{eqnarray}
&&\sum_{\mathbf{i}\in[ \partial R_n^{t_n} \setminus\partial
R_n^{u_n}]\cap\mathbb{Z}^d} \theta_n(
\mathbf{i})^2\nonumber
\\[-2pt]
&&\quad \leq C(d,\beta) \sum_{\mathbf{i}\in [ \partial R_n^{t_n}
\setminus\partial R_n^{u_n}]\cap\mathbb{Z}^d} \bigl\{
v_n(\mathbf{i})^{d-\beta}L\bigl(v_n(\mathbf{i})
\bigr) \bigr\}^2\nonumber
\\[-2pt]
\label{i-V3n}&& \quad \leq  C(d,\beta) \max \bigl\{ L(\mathbf{i})^2 \dvt  u_n
\leq\|\mathbf{i}\|\leq t_n \bigr\}\cdot \lambda_n^{3d-2\beta}
\int_{\partial R_0^{t_n/\lambda_n}} d(\mathbf{x}, \partial R_0)^{2d-2\beta}
\,\mathrm{d}\mathbf{x}
\\[-2pt]
&& \quad \leq C(d,\beta) \max \bigl\{ L(\mathbf{i})^2 \dvt  u_n
\leq\|\mathbf{i}\|\leq t_n \bigr\} \cdot \lambda_n^{3d-2\beta}
\int_0^{t_n/\lambda_n} t^{2d-2\beta}\, \mathrm{d}t\nonumber
\\[-2pt]
&& \quad = \mathrm{o} \bigl( \lambda_n^{3d-2\beta}L(\lambda_n)^2
\bigr).\nonumber
\end{eqnarray}
Hence, by (\ref{i-v3n}) and (\ref{i-V3n}),\vspace*{-3pt} it follows that
\[
\sum_{\mathbf{i}\in V_{3n}} \theta_n(
\mathbf{i})^2 = \mathrm{o} \bigl({\lambda}_{n}^{[3d -2\beta]}L({
\lambda}_{n})^2 \bigr).
\]
Next,\vspace*{-3pt} using arguments similar to those in the proof of Theorem~\ref{th32},
we have
\begin{eqnarray*}
\sum_{\mathbf{i}\in V_{1n}} \theta_n(
\mathbf{i})^2 &=& \sum_{\mathbf{i}\in V_{1n}} \biggl[\sum
_{\mathbf{j}\in R_n
-\mathbf{i}} \alpha(\mathbf{i}) \biggr]^2
\\
&=& \sum_{\mathbf{i}\in V_{1n}} \biggl[\sum
_{\mathbf{j}\in R_n^c
-\mathbf{i}} \alpha(\mathbf{i}) \biggr]^2\qquad \mbox{(as
$A=0$)}
\\
&=& {\lambda}_{n}^{2d}\gamma({\lambda}_{n})^2
\sum_{\mathbf{i}\in V_{1n}} \biggl[\int g_{{\lambda}_{n}}(\mathbf{x})
\mathbh{1} \bigl(\mathbf {x}\in  \bigl\llbracket  {R_0^c - {
\lambda}_{n}^{-1}{\mathbf{i}}}\bigr\rrbracket_n \bigr) \,\mathrm{d}\mathbf{x}
\biggr]^2,
\end{eqnarray*}
and\vspace*{-3pt} similarly,
\[
\sum_{\mathbf{i}\in V_{2n}} \theta_n(
\mathbf{i})^2 = {\lambda}_{n}^{2d}\gamma({
\lambda}_{n})^2 \sum_{\mathbf{i}\in V_{2n}}
\biggl[\int g_{{\lambda}_{n}}(\mathbf{x}) \mathbh{1} \bigl(\mathbf {x}\in
\bigl\llbracket {R_0 - {\lambda}_{n}^{-1}{\mathbf{i}}}
\bigr\rrbracket_n\bigr)\, \mathrm{d}\mathbf{x} \biggr]^2,
\]
where $\llbracket R_0 - {\lambda}_{n}^{-1}{\mathbf{i}}\rrbracket_n$ is as defined
in the
proof of Theorem~\ref{th32}.
Define the function $\check{G}_n(\cdot)$\vspace*{-3pt} by
%
\begin{equation}\label{chG}
\check{G}_n(\mathbf{x}) = \cases{
\ds \int
g_{{\lambda}_{n}}(\mathbf{y}) \mathbh{1} \bigl(\mathbf{y}\in \bigl\llbracket R_0^c - {\lambda}_{n}^{-1}{
\mathbf{i}}\bigr\rrbracket_n \bigr)\, \mathrm{d}\mathbf{y}, &\quad $\mbox{if } \mathbf{x}\in{
\lambda}_{n}^{-1}(\mathbf{i}+\mathcal {C}), \mathbf{i}\in
V_{1n}$,
\vspace*{3pt}\cr
\ds\int g_{{\lambda}_{n}}(\mathbf{y}) \mathbh{1} \bigl(\mathbf{y}\in  \bigl\llbracket R_0 - {\lambda}_{n}^{-1}{\mathbf{i}}
\bigr\rrbracket_n  \bigr) \,\mathrm{d}\mathbf{y} &\quad$\mbox{if }\mathbf{x}\in{\lambda}_{n}^{-1}(
\mathbf{i}+\mathcal {C}), \mathbf{i}\in V_{2n}$,
\vspace*{3pt}\cr
0, &\quad$\mbox{otherwise}$.}
\end{equation}
Then,
it follows\vspace*{-4pt} that
%
\begin{eqnarray}
\sigma_n^2 &= & {\lambda}_{n}^{3d}
\gamma({\lambda}_{n})^2 \int \check{G}_n(
\mathbf{x})^2 \,\mathrm{d}\mathbf{x}+ \sum_{\mathbf{i}\in V_{3n}} \bigl|
\theta_n(\mathbf{i})\bigr|^2
\nonumber
\\[-8pt]
\label{si-d-half}\\[-8pt]
\nonumber
&= & {\lambda}_{n}^{3d}\gamma({\lambda}_{n})^2
\int\check {G}_n(\mathbf{x})^2\, \mathrm{d}\mathbf{x}+ \mathrm{o} \bigl( {
\lambda}_{n}^{3d}\gamma({\lambda}_{n})^2
\bigr).
\end{eqnarray}
It now remains to show that
$\int\check{G}_n(\mathbf{x})^2\, \mathrm{d}\mathbf{x}\to
\int{G^\dagger_{\infty}}(\mathbf{x})^2\, \mathrm{d}\mathbf{x}$, or equivalently,
that $\int[\check{G}_n(\mathbf{x})-{G^\dagger_{\infty}}(\mathbf
{x})]^2\, \mathrm{d}\mathbf{x}\to
0$.
To that end, define
$
\Gamma_{n} = \bigcup\{{\lambda}_{n}^{-1}(\mathbf{i}+\mathcal{C})\dvt
\mathbf{i}\in V_{2n} \}$ and recall that $a_n=\sqrt{d}/\lambda_n$. Then,
repeating the arguments leading to
(\ref{diff1}),\vspace*{-2pt} we get
\begin{eqnarray*}
&&\sum_{\mathbf{i}\in V_{2n}} \int_{{\lambda}_{n}^{-1}(\mathbf{i}+\mathcal{C})} \biggl(
\int \bigl[ g_{{\lambda}_{n}}(\mathbf{x}) - {g}_{\infty}(\mathbf {x})\bigr]
\mathbh{1} \bigl(\mathbf{x}\in \bigl\llbracket R_0 - {\lambda}_{n}^{-1}{
\mathbf{i}}\bigr\rrbracket_n \bigr)\,\mathrm{d}\mathbf {x} \biggr)^2 \,\mathrm{d}\mathbf{z}
\\
&& \quad \leq \sum_{\mathbf{i}\in V_{2n} } \int_{{\lambda}_{n}^{-1}(\mathbf
{i}+\mathcal{C})}
\biggl( \int \bigl| g_{{\lambda}_{n}}^{(n)}(\mathbf{x}) - {g}_{\infty
}^{(n)}(
\mathbf{x}) \bigr| \mathbh{1} \bigl(\mathbf{x}\in R_0^{2a_n} -
\mathbf{z} \bigr)\,\mathrm{d}\mathbf{x} \biggr)^2 \,\mathrm{d}\mathbf{z}
\\
&&\quad= 
\int_{\Gamma_{n}} \biggl( \int \bigl|
g_{{\lambda}_{n}}^{(n)}(\mathbf{x}-\mathbf{z}) - {g}_{\infty}^{(n)}(
\mathbf{x}-\mathbf{z}) \bigr| \mathbh{1} \bigl(\mathbf{x}\in R_0^{2a_n}
\bigr)\,\mathrm{d}\mathbf{x} \biggr)^2 \,\mathrm{d}\mathbf{z}
\\
&& \quad= \int\!\!\!\int\!\!\!\int \bigl| g_{{\lambda}_{n}}^{(n)}(\mathbf{x}-\mathbf{z}) -
{g}_{\infty}^{(n)}(\mathbf{x}-\mathbf{z}) \bigr| \bigl| g_{{\lambda}_{n}}^{(n)}(
\mathbf{y}-\mathbf{z}) - {g}_{\infty
}^{(n)}(\mathbf{y}-\mathbf{z})
\bigr|
\\
&&\hspace*{23pt}\qquad{}\times\mathbh{1} \bigl(\mathbf{x}\in R_0^{2a_n}
\bigr) \mathbh{1} \bigl(\mathbf{y}\in R_0^{2a_n} \bigr)
\mathbh{1} (\mathbf{z}\in\Gamma_{n} )\, \mathrm{d}\mathbf{x} \,\mathrm{d}\mathbf{y}\, \mathrm{d}
\mathbf{z}
\\
&& \quad \leq \int\!\!\!\int \biggl(\int_{\Gamma_n} \bigl| g_{{\lambda}_{n}}^{(n)}(
\mathbf {x}-\mathbf{z}) - {g}_{\infty}^{(n)}(\mathbf{x}-\mathbf{z})
\bigr|^2 \,\mathrm{d}\mathbf{z} \biggr)^{1/2} \biggl(\int
_{\Gamma_n} \bigl| g_{{\lambda}_{n}}^{(n)}(\mathbf{y}-
\mathbf{z}) - {g}_{\infty
}^{(n)}(\mathbf{y}-\mathbf{z})
\bigr|^2 \,\mathrm{d}\mathbf{z} \biggr)^{1/2}
\\
&&\hspace*{16pt}\qquad{}\times \mathbh{1} \bigl(\mathbf{x}\in R_0^{2a_n}
\bigr) \mathbh{1} \bigl(\mathbf{y}\in R_0^{2a_n} \bigr)\, \mathrm{d}
\mathbf{x}\, \mathrm{d}\mathbf{y}
\\
&&\quad = \bigl\| \bigl(g_{{\lambda}_{n}}^{(n)} - {g}_{\infty}^{(n)}
\bigr) 
 \bigr\|_2^2 \biggl[\int\mathbh{1} \bigl(
\mathbf{x}\in R_0^{2a_n} \bigr)\,\mathrm{d}\mathbf{x}
\biggr]^2
\\
&&\quad=\mathrm{o}(1).
\end{eqnarray*}

By similar arguments,
\begin{eqnarray*}
&&\sum_{\mathbf{i}\in V_{1n}} \int_{{\lambda}_{n}^{-1}(\mathbf{i}+\mathcal{C})} \biggl(
\int \bigl[ g_{{\lambda}_{n}}(\mathbf{x}) - {g}_{\infty}(\mathbf {x})\bigr]
\mathbh{1} \bigl(\mathbf{x}\in \bigl\llbracket R_0^c - {\lambda}_{n}^{-1}{\mathbf{i}}\bigr\rrbracket_n \bigr)\,\mathrm{d}\mathbf{x}
\biggr)^2 \,\mathrm{d}\mathbf{z}
\\
&&\quad \leq  C(d) \cdot \biggl[ \int_{\{\|\mathbf{y}\|\geq{\eta}_{n}\}} \bigl| g_{{\lambda}_{n}} (
\mathbf{y}) - {g}_{\infty}(\mathbf{y}) \bigr| \,\mathrm{d}\mathbf{y} \biggr]^2 =
\mathrm{o}(1).
\end{eqnarray*}

Next, note that $\partial[R_0^c] = \partial R_0$.
By repeating the arguments in (\ref{diff2}), one can conclude that
\begin{eqnarray*}
&& \sum_{\mathbf{i}\in V_{1n}} \int_{{\lambda}_{n}^{-1}(\mathbf
{i}+\mathcal{C})}
\biggl( \int {g}_{\infty}(\mathbf{x}) \bigl[\mathbh{1} \bigl(\mathbf{x}\in
\bigl\llbracket  R_0^c - {\lambda}_{n}^{-1}{
\mathbf{i}}\bigr\rrbracket_n \bigr) - \mathbh{1} \bigl(\mathbf{x}\in
R_0^c - \mathbf{z} \bigr) \bigr]\, \mathrm{d}\mathbf{x}
\biggr)^2 \,\mathrm{d}\mathbf{z}
\\
&& \qquad {}+\sum_{\mathbf{i}\in V_{2n} } \int
_{{\lambda}_{n}^{-1}(\mathbf
{i}+\mathcal{C})} \biggl( \int {g}_{\infty}(\mathbf{x}) \bigl[
\mathbh{1} \bigl(\mathbf{x}\in \bigl\llbracket R_0 - {\lambda}_{n}^{-1}{
\mathbf{i}}\bigr\rrbracket_n \bigr) - \mathbh{1} (\mathbf{x}\in R_0 -
\mathbf{z} ) \bigr] \,\mathrm{d}\mathbf{x} \biggr)^2 \,\mathrm{d}\mathbf{z}
\\
&&\quad=\mathrm{o}(1).
\end{eqnarray*}

Finally, using the boundary condition (C.3) and the bounds on
$|{G^\dagger_{\infty}}(\mathbf{x})|$ for $\mathbf{x}\in(\partial
R_0)^{{\eta}_{n}}$ (from the
proof of ${G^\dagger_{\infty}}\in L^2(\mathbb{R}^d)$), one gets
$\sum_{\mathbf{i}\in V_{3n}}
\int_{{\lambda}_{n}^{-1}(\mathbf{i}+\mathcal{C})} cGi^2(\mathbf
{x})\,\mathrm{d}\mathbf{x}\leq C(d)
\int_0^{{\eta}_{n}} t^{2(d-\beta)}\,\mathrm{d}t =\mathrm{o}(1)$.
Hence,
it follows that
\[
\sigma_n^2= \bigl[\gamma({\lambda}_{n})
\bigr]^2 {\lambda}_{n}^{3d}\int
\bigl[{G^\dagger _{\infty}}(\mathbf{x})\bigr]^2\, \mathrm{d}
\mathbf{x}\bigl(1+\mathrm{o}(1)\bigr). 
\]
This completes the proof of Theorem~\ref{th41}.
\end{pf*}

\begin{pf*}{Proofs of claims in Example~\ref{ex42}}
Let $c_n= \lfloor{\log{\lambda}_{n}}\rfloor$, $n\geq1$. Also, let
$\|(x,y)'\|_\infty=\max\{|x|,  |y|\}$, $x,y\in\mathbb{R}$.
For a set $A\subset\mathbb{R}^2$ and $\delta>0$, write
$A^{\delta}_\infty= \{\mathbf{x}\in\mathbb{R}^2\dvt  \|\mathbf
{x}-\mathbf{y}\|_\infty
\leq\delta\}$ and $A^{-\delta}_\infty=
\{\mathbf{x}\in A \dvt \|\mathbf{x}-\mathbf{y}\|_\infty<\delta$ implies
$\mathbf{y}\in A\}$
for the $\delta$-enlargement and the $\delta$-interior
of a set $A$ in the $\|\cdot\|_\infty$-norm.
As before, set
$\theta_n(\mathbf{i}) = \sum_{\mathbf{j}+\mathbf{i}\in R_n}\alpha
(\mathbf{j})$.
Let $\mathcal{I}_n(i_0) \equiv
\{i\in\mathbb{Z}\dvt  -\frac{\lambda_n}{2} - i_0<i< \frac{\lambda
_n}{2} - i_0\}$,
$i_0\in\mathbb{Z}$.
For $\mathbf{i}=(i_0,j_0)\in R_{n,\infty}^{-c_n}$,
it is easy to verify that both $\mathcal{I}_n(i_0)$ and $\mathcal{I}_n(j_0)$
contain the set $\{i\in\mathbb{Z}\dvt  |i| < c_n\}$. Hence,
using the fact that $\sum_{\mathbf{i}\in\mathbb{Z}^2} \alpha
(\mathbf{i}) =0$, one gets
\begin{eqnarray*}
\bigl|\theta_n(\mathbf{i})\bigr| = \bigl|\theta_n(i_0,j_0)\bigr|&=&
\biggl| \sum_{i\in\mathcal{I}_n(i_0)\times\mathcal{I}_n(j_0)} \alpha(i,j) \biggr|
\\
&=& \biggl| \sum_{(i,j)\notin\mathcal{I}_n(i_0)\times\mathcal
{I}_n(j_0)} \alpha(i,j) \biggr|
\\
&=& \biggl| \sum_{(i,j)\notin\mathcal{I}_n(i_0)\times\mathcal
{I}_n(j_0)} b(i)b(j) \mathbh{1}(ij\neq0) \biggr|
\\
&\leq& \biggl[\sum_{|i|> {\lambda_n}/{2} - |i_0|} b(i) \biggr] \cdot \biggl[
\sum_{|i|> {\lambda_n}/{2} - |j_0|} b(i) \biggr].
\end{eqnarray*}
Hence, it follows that
%
\begin{eqnarray}
\sum_{\mathbf{i}\in R_{n,\infty}^{-c_n}} \theta_n(
\mathbf{i})^2 &\leq& \Biggl[\sum_{i_0= -\lambda_n/2+c_n}^{\lambda_n/2-c_n}
\biggl\{\sum_{|i|> {\lambda_n}/{2} - |i_0|} b(i) \biggr\}^2
\Biggr]^2\nonumber
\\
\label{E4-ub-1}&\leq& C(\beta,B) \Biggl[\sum_{j=0}^{[\lambda_n/2] -c_n} \biggl|
\frac{\lambda_n}{2} -j \biggr|^{2-2\beta} \Biggr]^2
\\
&\leq& C(\beta,B) \bigl(c_n^{3-2\beta} \bigr)^2\nonumber
%
\end{eqnarray}
for some (generic) constant $C(\beta,B)\in(0,\infty)$. By similar arguments,
it can be shown that
\[
\sum_{\mathbf{i}\notin R_{n,\infty}^{c_n}} \theta_n(\mathbf
{i})^2=\mathrm{o}(\lambda_n).
\]
Hence, it remains to determine the contribution of the boundary terms
to $\sigma_n^2$. For $ |i_0|\leq\lambda_n/2 -c_n$ and
$-c_n\leq k\equiv\lfloor{2^{-1}\lambda_n - j_0}\rfloor\leq c_n$
(this corresponds to
a part of the upper boundary line
of $R_n$), note that
\begin{eqnarray*}
\theta_n(i_0,j_0) &=& \sum
_{|i+i_0|< \lambda_n/2} \sum_{j=-{\lambda
}_{n}/2 -j_0}^k
\alpha(i,j)
\\
&=& \Biggl[ \sum_{i\in\mathbb{Z}}\sum
_{j=-{\lambda}_{n}/2 -j_0}^k \alpha(i,j) \Biggr]\bigl(1+\mathrm{o}(1)\bigr)
\\
&=& \mathrm{o}(1)+ \cases{
\ds -4B^2 + (2B) \Biggl[\sum
_{j=1}^k b(j) +B \Biggr], &\quad$\mbox{if }k > 0$,
\vspace*{3pt}\cr
-2B^2, &\quad$\mbox{if }k= 0$,
\vspace*{3pt}\cr
\ds (2B)\Biggl[\sum_{j=-\infty}^k b(j) \Biggr], &
\quad$\mbox{if }k<0$}
\end{eqnarray*}
uniformly in $(i_0,j_0)$. Note that by absolute summability of the
$\alpha(\mathbf{i})$,
\[
\sum_{|i_0+{\lambda}_{n}/2|\leq c_n}\sum_{|j_0-\lambda_n/2|\leq
c_n}
\theta_n(i_0,j_0)^2 = \mathrm{O}
\bigl(c_n^2\bigr).
\]

Hence, it follows that
\[
\sum_{| i_0|\leq\lambda_n/2 - c_n} \sum_{|j_0-\lambda_n/2|\leq c_n}
\theta_n(i_0,j_0)^2 = \bigl[ {
\lambda}_{n}\sigma_0^2/4 \bigr]\bigl(1+\mathrm{o}(1)
\bigr).
\]
Now using similar arguments
for the other three boundary arms, one gets the results of Example~\ref{ex42}.
\end{pf*}

\begin{pf*}{Proof of Theorem~\ref{th43}}
Using (\ref{e-cond}), one can show that
there exists ${\eta}_{n}\to0+$ such that
%
\begin{equation}\label{ee-cond}
\limsup_{n\to\infty}  \biggl|{\lambda}_{n}^{-(d-1)} \sum
_{\mathbf{i}\in[\partial R_n]^{{\eta}_{n}{\lambda}_{n}}\cap
\mathbb{Z}^d} \bigl|\theta_n(\mathbf{i})\bigr|^2
- \sigma_{\mathrm{EE}}^2 \biggr| =0.
\end{equation}
By (\ref{e-cond}) and the monotonicity of the sum of $\theta
_n(\mathbf{i})^2$
over increasing index sets, we may always
replace ${\eta}_{n}$ by a coarser sequence going to zero.
Hence, w.l.o.g. assume that
${\eta}_{n}\geq [\log(\lambda_n)]^{-1}$.
Next set $t_n = {\lambda}_{n}{\eta}_{n}$ and
(re-)define the sets $V_{kn},  k=1,2,3 $ in (\ref{v-s}) with
this choice of $t_n$. Further, write
$V_{21n} = V_{2n}\cap\{\mathbf{i}\in\mathbb{Z}^d\dvt  \|\mathbf{i}\|
\leq\lambda_n\sqrt{d}\}$
and $V_{22n} = V_{2n}\setminus V_{21n}$.

Next, note that uniformly in $\mathbf{i}\notin R_n^{t_n}$,
%
\begin{eqnarray}
\bigl|\theta_n(\mathbf{i})\bigr| &\leq& \sum_{\mathbf{j}\in R_n-\mathbf{i}}
\bigl|\alpha(\mathbf{j})\bigr|\nonumber
\\
\label{tb-1}&\leq &  \sum_{\mathbf{j}\dvt  \|\mathbf{j}\| \geq d(\partial R_n,\mathbf{i}),
 \mathbf{j}\in[R_n-\mathbf{i}]\cap\mathbb{Z}^d} \gamma\bigl(\|\mathbf {j}\|\bigr)
\\
& \leq &   C(d,\beta) \min \bigl\{ d(\partial R_n,\mathbf{i})^{d-\beta}
L\bigl(d(\partial R_n,\mathbf{i})\bigr),  \lambda_n^d
d(\partial R_n,\mathbf{i})^{-\beta} L^*_{n}(
\mathbf{i}) \bigr\},\nonumber
\end{eqnarray}
where $ L^*_{n}(\mathbf{i}) \equiv\max\{ L(\|\mathbf{j}\|) \dvt
\mathbf{j}\in[R_n-\mathbf{i}]\cap\mathbb{Z}^d\}$.
And, using the fact that $A=0$, one can similarly show that
uniformly in $\mathbf{i}\in R_n^{-t_n}$,
\[
\bigl|\theta_n(\mathbf{i})\bigr| = \biggl| - \sum_{\mathbf{j}\in R_n^c-\mathbf
{i}}
\alpha(\mathbf{j}) \biggr| \leq C(d,\beta) d\bigl(\partial R_n^c,
\mathbf{i}\bigr)^{d-\beta} L\bigl(d\bigl(\partial R_n^c,
\mathbf{i}\bigr)\bigr).
\]
Hence, it follows that
%
\begin{eqnarray}
&&\sum_{\mathbf{i}\in V_{1n}}\theta_n(
\mathbf{i})^2 +\sum_{\mathbf{i}\in V_{21n}}
\theta_n(\mathbf{i})^2\nonumber
\\
\label{v12} && \quad\leq C(d) \lambda_n^d \cdot\max\bigl\{\theta_n(\mathbf{i})^2 \dvt  \mathbf{i}\in V_{1n}
\cup V_{21n}\bigr\}
\\
&& \quad\leq C(d,\beta) \lambda_n^d t_n^{2(d-\beta)}
\max\bigl\{L^2(t)\dvt  t_n\leq t\leq d\lambda_n
\bigr\} = \mathrm{o}\bigl(\lambda_n^{d-1}\bigr).\nonumber
\end{eqnarray}
%
Next, note that for all $\mathbf{i}\in V_{22n}$,
\[
\|\mathbf{i}\|/2\leq\|\mathbf{i}\|-\lambda_n\sqrt{d}/2 \leq d(
\partial R_n, \mathbf{i}) \leq\|\mathbf{i}\|
\]
and
\[
\max\bigl\{\|\mathbf{j}\|\dvt  \mathbf{j}\in[R_n-
\mathbf{i}]\cap\mathbb {Z}^d\bigr\} \subset\bigl[\|\mathbf{i}\|/2, 3\|
\mathbf{i}\|/2\bigr\}.
\]
Hence, by (\ref{tb-1})
%
\begin{eqnarray}
\sum_{\mathbf{i}\in V_{22n}}\theta_n(
\mathbf{i})^2 &\leq& C(d) \lambda_n^{2d} \sum
_{\|\mathbf{i}\|>\sqrt{d}\lambda
_n} \|\mathbf{i}\|^{-2\beta} L^2\bigl(
\|\mathbf{i}\|\bigr)
\nonumber
\\[-8pt]
\label{v22}\\[-8pt]
\nonumber
&\leq& C(d,\beta) \lambda_n^{3d-2\beta}L^2(
\lambda_n) = \mathrm{o}\bigl(\lambda_n^{d-1}\bigr)
\end{eqnarray}
for $\beta>d+1/2$.
Hence, from (\ref{e-cond}), (\ref{v12}) and (\ref{v22}), it follows that
%
\begin{equation}\label{ee-th}
\sigma_n^2 = \sum_{k=1}^3
\sum_{\mathbf{i}\in V_{kn}} \theta _n(
\mathbf{i})^2 
= {\lambda}_{n}^{d-1}
\sigma^2_{\mathrm{EE}}\bigl(1+\mathrm{o}(1)\bigr).
\end{equation}
%
By Corollary~\ref{co62}(i), Theorem~\ref{th43} follows.
\end{pf*}

\begin{pf*}{Proof of Theorem~\ref{thm44}}
Theorem~\ref{thm44} follows from
Corollary~\ref{co62}(i), by comparing the orders of the
terms $\sum_{\mathbf{i}\in V_{kn}} \theta_n(\mathbf{i})^2$, $k=1,2,3$.
Specifically, for part (i), we use the arguments in the proof of
Theorem~\ref{th43},
and the bounds from (\ref{ee-cond}), (\ref{v12}) and (\ref{v22})
to conclude that (\ref{ee-th}) holds. For part (ii), note that
for $\beta=d+1/2$, the conditions and the arguments in the proof of
Theorem~\ref{th41} no longer
ensure that $\int_{\|\mathbf{x}\|\leq{\eta}_{n}} {G^\dagger
_{\infty}}(\mathbf{x})^2 \,\mathrm{d}\mathbf{x}= \mathrm{O}(1)$, which is why we need to
make the assumption that ${G^\dagger_{\infty}}\in L^2$. However,
under (C.1) and (C.2),
the arguments in the proof of Theorem~\ref{th41}
leading to the convergence of $\int_{ \|\mathbf{x}\| \geq{\eta
}_{n}} |\check{G}_n(\mathbf{x}) -{G^\dagger_{\infty}}(\mathbf
{x})|\,\mathrm{d}\mathbf{x}$
to zero still holds. Both parts of (ii) now follow by
deriving the limits of the terms $\sum_{\mathbf{i}\in V_{kn}} \theta
_n(\mathbf{i})^2$,
using (\ref{ee-cond}) for $k=3$ and using the steps from the proof of
Theorem~\ref{th41}
for $k=1,2$. We omit the routine details.
\end{pf*}

\subsection{Proof of the results from Section~\texorpdfstring{\protect\ref{sec5}}{5}}\label{sec65}

\begin{pf*}{Proof of Theorem~\ref{th51}}
In view of Corollary~\ref{co62}(ii),
it is enough to show that
%
\begin{equation}\label{si-srd}
\sigma_n^2 = N_n A^2 \bigl(1 +
\mathrm{o}(1)\bigr).
\end{equation}
Let $c_n$ be a sequence of positive real numbers satisfying
%
\begin{equation}\label{den}
c_n^{-1} + {\lambda}_{n}^{-1}c_n
= \mathrm{o}(1).
\end{equation}
Also, let
%
\begin{eqnarray}
U_{1n} &=&\bigl\{\mathbf{i}\in\mathbb{Z}^d \dvt  \bigl\| {
\lambda}_{n}^{-1}\mathbf {i}-\mathbf{x}\bigr\|\leq
c_n/{\lambda}_{n} \mbox{ for some } \mathbf{x}\in\partial
R_0\bigr\},
\nonumber\\
U_{2n} &=&\bigl\{\mathbf{i}\in\mathbb{Z}^d\dvt  B(
\mathbf{i};c_n)\subset R_n\bigr\},
\nonumber
\\[-8pt]
\label{U14}\\[-8pt]
\nonumber
U_{3n} &=& \bigl\{\mathbf{i}\in[- 2d{\lambda}_{n},2d{
\lambda}_{n}]^d \dvt  \mathbf{i}\notin[U_{1n}\cup
U_{2n}]\bigr\}\quad \mbox{and}
\\
U_{4n} & =& \mathbb{Z}^d\setminus U_{3n}.
\nonumber
\end{eqnarray}
%
Then $\sigma_n^2$ can be written as
%
\begin{equation}\label{srd-i14}
\sigma_n^2 = \sum_{i=1}^4
\sum_{\mathbf{i}\in U_{in}} \theta _n(
\mathbf{i})^2 \equiv I_{1n}+ I_{2n}+I_{3n}+I_{4n},\qquad
\mbox{say}.
\end{equation}
Note that by (\ref{den}), the boundary condition on $R_0$
(that $\nu(\partial R_0)=0$) and the absolute
summability of $\alpha(\mathbf{i})$'s,
\begin{eqnarray*}
I_{1n} &\leq& \biggl[\sum_{\mathbf{j}\in\mathbb{Z}^d}\bigl|\alpha (
\mathbf{j})\bigr| \biggr]^2 \bigl| \bigl\{\mathbf{i}\dvt  {\lambda}_{n}^{-1}
\mathbf{i}\in [\partial R_0]^{c_n/{\lambda}_{n}} \bigr\} \bigr|
\\
&=& {\lambda}_{n}^d \cdot \mathrm{O} \bigl(\operatorname{vol.}
\bigl( [\partial R_0]^{c_n/{\lambda}_{n}} \bigr) \bigr)
\\
& =& \mathrm{o}\bigl({\lambda}_{n}^d\bigr).
\end{eqnarray*}
Next, consider $I_{2n}$. By definition of $U_{2n}$, $B(\mathbf{0};c_n)
\subset
R_n - \mathbf{i}$ for all $\mathbf{i}\in U_{2n}$. Hence,
it follows that
%
\begin{equation}\label{i2}
\sup\bigl\{\bigl|\theta_n(\mathbf{i}) - A\bigr| \dvt  \mathbf{i}\in
U_{2n}\bigr\} \leq \sum_{\mathbf{j}\in\mathbb{Z}^d \dvt  \|\mathbf{j}\|\geq c_n} \bigl|\alpha (
\mathbf{j})\bigr| = \mathrm{o}(1).
\end{equation}
Next, note that for any $\mathbf{i}\in U_{3n}$,
$B(\mathbf{i};c_n/2)$ is contained in the set $R_n^c$, and hence,
$\varnothing=[R_n -\mathbf{i}]\cap B(\mathbf{0}; c_n/2) = [R_n\cap
B(\mathbf{i}; c_n/2)]
-\mathbf{i}$.
Hence, it follows that
\[
I_{3n} \leq C(d) \lambda_n^d \sup \bigl\{
\theta_n(\mathbf{i})^2\dvt  \mathbf {i}\in U_{3n}
\bigr\} \leq C(d) \lambda_n^d \sum
_{\mathbf{j}\in\mathbb{Z}^d\dvt  2 \|\mathbf{j}\|\geq c_n} \bigl|\alpha(\mathbf{j})\bigr| = \mathrm{o}\bigl({\lambda}_{n}^d
\bigr).
\]
Finally, by condition (C.1) and the definition of $U_{4n}$,
\begin{eqnarray*}
&&\sum_{\mathbf{i}\in U_{4n}} \theta_n(
\mathbf{i})^2
\\
&&\quad \leq \sum_{\mathbf{i}\in U_{4n}} \biggl[ \biggl(\sum
_{\mathbf{j}\in
[R_n-\mathbf{i}]\cap\mathbb{Z}^d} \alpha(\mathbf{j})^2 \biggr) \times
N_n \biggr]
\\
&& \quad\leq N_n \sum_{\mathbf{i}\in U_{4n}} \sum
_{\mathbf{j}\in\mathbb{Z}^d\dvt  \|
\mathbf{j}\|_1\geq
\|\mathbf{i}\|_1 -d{\lambda}_{n}} \alpha(\mathbf{j})^2,\qquad  \Bigl(\mbox{since}
\sup_{\mathbf{x}\in R_n}\|\mathbf{x}\|_1 \leq d{
\lambda}_{n}/2\Bigr)
\\
&&\quad\leq  N_n \sum_{\mathbf{i}\in U_{4n}} \sum
_{k\geq \|\mathbf{i}\|_1 -
d{\lambda}_{n}} \bigl|\bigl\{\mathbf{j}\in\mathbb{Z}^d\dvt \|
\mathbf{j}\|_1 = k\bigr\} \bigr| \cdot\sup_{ \mathbf{j}\in\mathbb{Z}^d\dvt \|\mathbf{j}\|_1 = k} \alpha(
\mathbf{j})^2
\\
&& \quad\leq C(d)
N_n \sum_{\mathbf{i}\in U_{4n}} \sum
_{k\geq\|\mathbf{i}\|_1 -
d{\lambda}_{n}} k^{d-1} \sup\bigl\{ \gamma(t)^2 \dvt
t\in[k/\sqrt{d}, k]\bigr\}
\\
&& \quad\leq C(d) N_n\sum_{\mathbf{i}\in U_{4n}} \int
_{\|\mathbf{i}\|_1
- d{\lambda}_{n}}^\infty t^{d-1-2\beta}L(t)^2\,\mathrm{d}t
\\
&&\quad\leq C(d) N_n\int_{2d{\lambda}_{n}}^\infty
u^{d-1} \int_{u - d{\lambda}_{n}}^\infty
t^{d-1-2\beta}L(t)^2\,\mathrm{d}t \,\mathrm{d}u
\\
&&\quad\leq C(d) N_n {\lambda}_{n}^{2d-2\beta}L({
\lambda}_{n})^2 = \mathrm{o}(N_n),
\end{eqnarray*}
since ${\lambda}_{n}^{d-\beta}L({\lambda}_{n}) =\mathrm{o}(1)$ for all $\beta
> d$ and also for $\beta=d$
by the integrability
condition $\int_1^\infty\gamma(t)\,\mathrm{d}t <\infty$.
Since $|U_{2n}| = |\mathcal{D}_n| - \mathrm{O}(|U_{1n}|) = N_n - \mathrm{o}({\lambda
}_{n}^d) = N_n(1+\mathrm{o}(1))$,
(\ref{si-srd}) follows from (\ref{srd-i14}), (\ref{i2}) and the
bounds for $I_{1n}$, $I_{3n}$ and $I_{4n}$ above. This completes the
proof of the theorem.
\end{pf*}

\section*{Acknowledgements}
Research partially supported by US NSF Grants DMS 1007703
and DMS 1310068
and
ESRC Grant ES/J007242/1.
The authors thank the Associate Editor
and a referee for a number of insightful and constructive
comments that
significantly improved an earlier draft of the paper.





%
%

\printhistory

\begin{thebibliography}{38}

\bibitem{AthLah06}
\begin{bbook}[mr]
\bauthor{\bsnm{Athreya},~\bfnm{Krishna~B.}\binits{K.B.}} \AND
\bauthor{\bsnm{Lahiri},~\bfnm{Soumendra~N.}\binits{S.N.}}
(\byear{2006}).
\btitle{Measure Theory and Probability Theory}.
\blocation{New York}:
\bpublisher{Springer}.
\bid{mr={2247694}}
\end{bbook}
%

\bptok{imsref}%
\endbibitem

\bibitem{BerGhoSch09}
\begin{barticle}[mr]
\bauthor{\bsnm{Beran},~\bfnm{Jan}\binits{J.}},
\bauthor{\bsnm{Ghosh},~\bfnm{Sucharita}\binits{S.}} \AND
\bauthor{\bsnm{Schell},~\bfnm{Dieter}\binits{D.}}
(\byear{2009}).
\btitle{On least squares estimation for long-memory lattice processes}.
\bjournal{J. Multivariate Anal.}
\bvolume{100}
\bpages{2178--2194}.
\bid{doi={10.1016/j.jmva.2009.04.007}, issn={0047-259X}, mr={2560362}}
\end{barticle}
%

\bptok{imsref}%
\endbibitem

\bibitem{bert06}
\begin{bbook}[mr]
\beditor{\bsnm{Bertail},~\bfnm{P.}\binits{P.}},
\beditor{\bsnm{Doukhan},~\bfnm{P.}\binits{P.}} \AND
\beditor{\bsnm{Soulier},~\bfnm{P.}\binits{P.}}
(\byear{2006}).
\btitle{Dependence in Probability and Statistics}.
\bseries{Lecture Notes in Statistics}
\bvolume{187}.
\blocation{New York}:
\bpublisher{Springer}.
\bid{doi={10.1007/0-387-36062-X}, mr={2269087}}
\end{bbook}
%

\bptok{imsref}%
\endbibitem

\bibitem{Bil95}
\begin{bbook}[mr]
\bauthor{\bsnm{Billingsley},~\bfnm{Patrick}\binits{P.}}
(\byear{1995}).
\btitle{Probability and Measure},
\bedition{3rd} ed.
\blocation{New York}:
\bpublisher{Wiley}.
\bid{mr={1324786}}
\end{bbook}
%

\bptok{imsref}%
\endbibitem

\bibitem{Boietal05}
\begin{barticle}[mr]
\bauthor{\bsnm{Boissy},~\bfnm{Y.}\binits{Y.}},
\bauthor{\bsnm{Bhattacharyya},~\bfnm{B.~B.}\binits{B.B.}},
\bauthor{\bsnm{Li},~\bfnm{X.}\binits{X.}} \AND
\bauthor{\bsnm{Richardson},~\bfnm{G.~D.}\binits{G.D.}}
(\byear{2005}).
\btitle{Parameter estimates for fractional autoregressive spatial processes}.
\bjournal{Ann. Statist.}
\bvolume{33}
\bpages{2553--2567}.
\bid{doi={10.1214/009053605000000589}, issn={0090-5364}, mr={2253095}}
\end{barticle}
%

\bptok{imsref}%
\endbibitem

\bibitem{Bol82}
\begin{barticle}[mr]
\bauthor{\bsnm{Bolthausen},~\bfnm{E.}\binits{E.}}
(\byear{1982}).
\btitle{On the central limit theorem for stationary mixing random fields}.
\bjournal{Ann. Probab.}
\bvolume{10}
\bpages{1047--1050}.
\bid{issn={0091-1798}, mr={0672305}}
\end{barticle}
%

\bptok{imsref}%
\endbibitem

\bibitem{CarMejMor85}
\begin{barticle}[auto:parserefs-M02]
\bauthor{\bsnm{Carlos-Davila},~\bfnm{E.}\binits{E.}},
\bauthor{\bsnm{Mejia-Lira},~\bfnm{F.}\binits{F.}} \AND
\bauthor{\bsnm{Moran-Lopez},~\bfnm{J.~L.}\binits{J.L.}}
(\byear{1985}).
\btitle{Ferromagnetism and spatial long range order in binary alloys: Systems with face-centered cubic structure}.
\bjournal{J. Phys. C: Solid State Phys.}
\bvolume{18}
\bpages{1217--1224}.
\end{barticle}
%

\bptok{imsref}%
\endbibitem

\bibitem{Cre93}
\begin{bbook}[mr]
\bauthor{\bsnm{Cressie},~\bfnm{Noel~A.~C.}\binits{N.A.C.}}
(\byear{1993}).
\btitle{Statistics for Spatial Data}.
\blocation{New York}:
\bpublisher{Wiley}.
\bid{mr={1239641}}
\end{bbook}
%

\bptok{imsref}%
\endbibitem

\bibitem{DobMaj79}
\begin{barticle}[mr]
\bauthor{\bsnm{Dobrushin},~\bfnm{R.~L.}\binits{R.L.}} \AND
\bauthor{\bsnm{Major},~\bfnm{P.}\binits{P.}}
(\byear{1979}).
\btitle{Non-central limit theorems for nonlinear functionals of {G}aussian fields}.
\bjournal{Z. Wahrsch. Verw. Gebiete}
\bvolume{50}
\bpages{27--52}.
\bid{doi={10.1007/BF00535673}, issn={0044-3719}, mr={0550122}}
\end{barticle}
%

\bptok{imsref}%
\endbibitem

\bibitem{Dou94}
\begin{bbook}[mr]
\bauthor{\bsnm{Doukhan},~\bfnm{Paul}\binits{P.}}
(\byear{1994}).
\btitle{Mixing: Properties and Examples}.
\bseries{Lecture Notes in Statistics}
\bvolume{85}.
\blocation{New York}:
\bpublisher{Springer}.
\bid{doi={10.1007/978-1-4612-2642-0}, mr={1312160}}
\end{bbook}
%

\bptok{imsref}%
\endbibitem

\bibitem{ElMVolWu13}
\begin{barticle}[mr]
\bauthor{\bsnm{El Machkouri},~\bfnm{Mohamed}\binits{M.}},
\bauthor{\bsnm{Voln{\'y}},~\bfnm{Dalibor}\binits{D.}} \AND
\bauthor{\bsnm{Wu},~\bfnm{Wei~Biao}\binits{W.B.}}
(\byear{2013}).
\btitle{A central limit theorem for stationary random fields}.
\bjournal{Stochastic Process. Appl.}
\bvolume{123}
\bpages{1--14}.
\bid{doi={10.1016/j.spa.2012.08.014}, issn={0304-4149}, mr={2988107}}
\end{barticle}
%

\bptok{imsref}%
\endbibitem

\bibitem{Fai38}
\begin{barticle}[auto:parserefs-M02]
\bauthor{\bsnm{Fairfield Smith},~\bfnm{H.}\binits{H.}}
(\byear{1938}).
\btitle{An empirical law describing heterogeneity in the yields of agricultural crops}.
\bjournal{J. Agricultural Sci.}
\bvolume{28}
\bpages{1--23}.
\end{barticle}
%

\bptok{imsref}%
\endbibitem

\bibitem{Gne00}
\begin{barticle}[mr]
\bauthor{\bsnm{Gneiting},~\bfnm{Tilmann}\binits{T.}}
(\byear{2000}).
\btitle{Power-law correlations, related models for long-range dependence and their simulation}.
\bjournal{J. Appl. Probab.}
\bvolume{37}
\bpages{1104--1109}.
\bid{issn={0021-9002}, mr={1808873}}
\end{barticle}
%

\bptok{imsref}%
\endbibitem

\bibitem{GomHaz70}
\begin{btechreport}[auto:parserefs-M02]
\bauthor{\bsnm{Gomez},~\bfnm{M.}\binits{M.}} \AND
\bauthor{\bsnm{Hazen},~\bfnm{K.}\binits{K.}}
(\byear{1970}).
\btitle{Evaluating sulfur and ash distribution in coal seams by statistical response surface regression
analysis}.
\btype{Report RI 7377, U.S. Bureau of Mines}.
\end{btechreport}
%

\bptok{imsref}%
\endbibitem

\bibitem{Guy95}
\begin{bbook}[mr]
\bauthor{\bsnm{Guyon},~\bfnm{Xavier}\binits{X.}}
(\byear{1995}).
\btitle{Random Fields on a Network: Modeling, Statistics, and
Applications}.
\blocation{New York}:
\bpublisher{Springer}.
\bid{mr={1344683}}
\end{bbook}
%

\bptok{imsref}%
\endbibitem

\bibitem{Han79}
\begin{barticle}[mr]
\bauthor{\bsnm{Hannan},~\bfnm{E.~J.}\binits{E.J.}}
(\byear{1979}).
\btitle{The central limit theorem for time series regression}.
\bjournal{Stochastic Process. Appl.}
\bvolume{9}
\bpages{281--289}.
\bid{doi={10.1016/0304-4149(79)90050-4}, issn={0304-4149}, mr={0562049}}
\end{barticle}
%

\bptok{imsref}%
\endbibitem

\bibitem{IbrLin71}
\begin{bbook}[mr]
\bauthor{\bsnm{Ibragimov},~\bfnm{I.~A.}\binits{I.A.}} \AND
\bauthor{\bsnm{Linnik},~\bfnm{Yu.~V.}\binits{Yu.V.}}
(\byear{1971}).
\btitle{Independent and Stationary Sequences of Random Variables}.
\blocation{Groningen}:
\bpublisher{Wolters-Noordhoff Publishing}.
\bid{mr={0322926}}
\end{bbook}
%

\bptok{imsref}%
\endbibitem

\bibitem{IvaLeo89}
\begin{bbook}[mr]
\bauthor{\bsnm{Ivanov},~\bfnm{A.~V.}\binits{A.V.}} \AND
\bauthor{\bsnm{Leonenko},~\bfnm{N.~N.}\binits{N.N.}}
(\byear{1989}).
\btitle{Statistical Analysis of Random Fields}.
\bseries{Mathematics and Its Applications (Soviet Series)}
\bvolume{28}.
\blocation{Dordrecht}:
\bpublisher{Kluwer Academic}.
\bid{doi={10.1007/978-94-009-1183-3}, mr={1009786}}
\end{bbook}
%

\bptok{imsref}%
\endbibitem

\bibitem{KasLap84}
\begin{barticle}[pbm]
\bauthor{\bsnm{Kashyap},~\bfnm{R.~L.}\binits{R.L.}} \AND
\bauthor{\bsnm{Lapsa},~\bfnm{P.~M.}\binits{P.M.}}
(\byear{1984}).
\btitle{Synthesis and estimation of random fields using long-correlation models}.
\bjournal{IEEE Trans. Pattern. Anal. Mach. Intell.}
\bvolume{6}
\bpages{800--809}.
\bid{issn={0162-8828}, pmid={22499661}}
\end{barticle}
%

\bptok{imsref}%
\endbibitem

\bibitem{Lah99}
\begin{barticle}[mr]
\bauthor{\bsnm{Lahiri},~\bfnm{S.~N.}\binits{S.N.}}
(\byear{1999}).
\btitle{Asymptotic distribution of the empirical spatial cumulative distribution function predictor and prediction bands based on a subsampling method}.
\bjournal{Probab. Theory Related Fields}
\bvolume{114}
\bpages{55--84}.
\bid{doi={10.1007/s004400050221}, issn={0178-8051}, mr={1697139}}
\end{barticle}
%

\bptok{imsref}%
\endbibitem

\bibitem{Lah03}
\begin{barticle}[mr]
\bauthor{\bsnm{Lahiri},~\bfnm{S.N.}\binits{S.N.}}
(\byear{2003}).
\btitle{Central limit theorems for weighted sums of a spatial process under a class of stochastic and fixed designs}.
\bjournal{Sankhy\=a}
\bvolume{65}
\bpages{356--388}.
\bid{issn={0972-7671}, mr={2028905}}
\bptnote{check year}%
\end{barticle}
%

\bptok{imsref}%
\endbibitem

\bibitem{Lahetal99}
\begin{barticle}[mr]
\bauthor{\bsnm{Lahiri},~\bfnm{Soumendra~N.}\binits{S.N.}},
\bauthor{\bsnm{Kaiser},~\bfnm{Mark~S.}\binits{M.S.}},
\bauthor{\bsnm{Cressie},~\bfnm{Noel}\binits{N.}} \AND
\bauthor{\bsnm{Hsu},~\bfnm{Nan-Jung}\binits{N.-J.}}
(\byear{1999}).
\btitle{Prediction of spatial cumulative distribution functions using subsampling}.
\bjournal{J. Amer. Statist. Assoc.}
\bvolume{94}
\bpages{86--110}.
\bid{doi={10.2307/2669680}, issn={0162-1459}, mr={1689216}}
\end{barticle}
%

\bptok{imsref}%
\endbibitem

\bibitem{Lav06}
\begin{bincollection}[mr]
\bauthor{\bsnm{Lavancier},~\bfnm{Fr{\'e}d{\'e}ric}\binits{F.}}
(\byear{2006}).
\btitle{Long memory random fields}.
In \bbooktitle{Dependence in Probability and Statistics}
(\beditor{\bfnm{P.}\binits{P.}~\bsnm{Bertail}},
\beditor{\bfnm{P.}\binits{P.}~\bsnm{Doukhan}} \AND
\beditor{\bfnm{P.}\binits{P.}~\bsnm{Soulier}}, eds.).
\bseries{Lecture Notes in Statist.}
\bvolume{187}
\bpages{195--220}.
\blocation{New York}:
\bpublisher{Springer}.
\bid{doi={10.1007/0-387-36062-X_9}, mr={2283256}}
\end{bincollection}
%

\bptok{imsref}%
\endbibitem

\bibitem{Lav07}
\begin{barticle}[mr]
\bauthor{\bsnm{Lavancier},~\bfnm{Fr{\'e}d{\'e}ric}\binits{F.}}
(\byear{2007}).
\btitle{Invariance principles for non-isotropic long memory random fields}.
\bjournal{Stat. Inference Stoch. Process.}
\bvolume{10}
\bpages{255--282}.
\bid{doi={10.1007/s11203-006-9001-9}, issn={1387-0874}, mr={2321311}}
\end{barticle}
%

\bptok{imsref}%
\endbibitem

\bibitem{LeoTau13}
\begin{barticle}[auto:parserefs-M02]
\bauthor{\bsnm{Leonenko},~\bfnm{N.}\binits{N.}} \AND
\bauthor{\bsnm{Taufer},~\bfnm{E.}\binits{E.}}
(\byear{2013}).
\btitle{Disaggregation of spatial autoregressive processes}.
\bjournal{Spatial Statistics}
\bvolume{3}
\bpages{1--20}.
\end{barticle}
%

\bptok{imsref}%
\endbibitem

\bibitem{MerHal11}
\begin{barticle}[auto:parserefs-M02]
\bauthor{\bsnm{Mercer},~\bfnm{W.~B.}\binits{W.B.}} \AND
\bauthor{\bsnm{Hall},~\bfnm{A.~D.}\binits{A.D.}}
(\byear{1911}).
\btitle{The experimental error of field trials}.
\bjournal{J. Agricultural Sci.}
\bvolume{4}
\bpages{107--132}.
\end{barticle}
%

\bptok{imsref}%
\endbibitem

\bibitem{Peretal08}
\begin{barticle}[auto:parserefs-M02]
\bauthor{\bsnm{Percival},~\bfnm{D.~B.}\binits{D.B.}},
\bauthor{\bsnm{Rothrock},~\bfnm{D.~A.}\binits{D.A.}},
\bauthor{\bsnm{Thorndike},~\bfnm{A.~S.}\binits{A.S.}} \AND
\bauthor{\bsnm{Gneiting},~\bfnm{T.}\binits{T.}}
(\byear{2008}).
\btitle{The variance of mean sea-ice thickness: Effect of long range dependence}.
\bjournal{J. Geophys. Research: Oceans}
\bvolume{113}.
\bid{doi={10.1029/2007JC004391}}
\end{barticle}
%

\bptok{imsref}%
\endbibitem

\bibitem{Res87}
\begin{bbook}[mr]
\bauthor{\bsnm{Resnick},~\bfnm{Sidney~I.}\binits{S.I.}}
(\byear{1987}).
\btitle{Extreme Values, Regular Variation, and Point Processes}.
\bseries{Applied Probability. A~Series of the Applied Probability Trust}
\bvolume{4}.
\blocation{New York}:
\bpublisher{Springer}.
\bid{doi={10.1007/978-0-387-75953-1}, mr={0900810}}
\end{bbook}
%

\bptok{imsref}%
\endbibitem

\bibitem{Rigetal87}
\begin{bbook}[auto:parserefs-M02]
\bauthor{\bsnm{Riggan},~\bfnm{W.~B.}\binits{W.B.}},
\bauthor{\bsnm{Creason},~\bfnm{J.~P.}\binits{J.P.}},
\bauthor{\bsnm{Nelson},~\bfnm{W.~C.}\binits{W.C.}},
\bauthor{\bsnm{Manton},~\bfnm{K.~G.}\binits{K.G.}},
\bauthor{\bsnm{Woodbury},~\bfnm{M.~A.}\binits{M.A.}},
\bauthor{\bsnm{Stallard},~\bfnm{E.}\binits{E.}},
\bauthor{\bsnm{Pellom},~\bfnm{A.~C.}\binits{A.C.}} \AND
\bauthor{\bsnm{Beaubier},~\bfnm{J.}\binits{J.}}
(\byear{1987}).
\btitle{U.S. Cancer Mortality Rates and Trends, 1950--1979 (Vol. IV: Maps)}.
\blocation{Washington, DC}:
\bpublisher{U.S. Environmental Protection Agency}.
\end{bbook}
%

\bptok{imsref}%
\endbibitem

\bibitem{Rob97}
\begin{barticle}[mr]
\bauthor{\bsnm{Robinson},~\bfnm{P.~M.}\binits{P.M.}}
(\byear{1997}).
\btitle{Large-sample inference for nonparametric regression with dependent errors}.
\bjournal{Ann. Statist.}
\bvolume{25}
\bpages{2054--2083}.
\bid{doi={10.1214/aos/1069362387}, issn={0090-5364}, mr={1474083}}
\end{barticle}
%

\bptok{imsref}%
\endbibitem

\bibitem{Ros56}
\begin{barticle}[mr]
\bauthor{\bsnm{Rosenblatt},~\bfnm{M.}\binits{M.}}
(\byear{1956}).
\btitle{A central limit theorem and a strong mixing condition}.
\bjournal{Proc. Natl. Acad. Sci. USA}
\bvolume{42}
\bpages{43--47}.
\bid{issn={0027-8424}, mr={0074711}}
\end{barticle}
%

\bptok{imsref}%
\endbibitem

\bibitem{Ros61}
\begin{bincollection}[mr]
\bauthor{\bsnm{Rosenblatt},~\bfnm{M.}\binits{M.}}
(\byear{1961}).
\btitle{Independence and dependence}.
In \bbooktitle{Proc. 4th {B}erkeley {S}ympos. {M}ath. {S}tatist. and {P}rob., {V}ol. {II}}
\bpages{431--443}.
\blocation{Berkeley, CA}:
\bpublisher{Univ. California Press}.
\bid{mr={0133863}}
\end{bincollection}
%

\bptok{imsref}%
\endbibitem

\bibitem{SheCar94}
\begin{barticle}[mr]
\bauthor{\bsnm{Sherman},~\bfnm{Michael}\binits{M.}} \AND
\bauthor{\bsnm{Carlstein},~\bfnm{Edward}\binits{E.}}
(\byear{1994}).
\btitle{Nonparametric estimation of the moments of a general statistic computed from spatial data}.
\bjournal{J. Amer. Statist. Assoc.}
\bvolume{89}
\bpages{496--500}.
\bid{issn={0162-1459}, mr={1294075}}
\end{barticle}
%

\bptok{imsref}%
\endbibitem

\bibitem{Sur82}
\begin{barticle}[mr]
\bauthor{\bsnm{Surgailis},~\bfnm{D.}\binits{D.}}
(\byear{1982}).
\btitle{Domains of attraction of self-similar multiple integrals}.
\bjournal{Litovsk. Mat. Sb.}
\bvolume{22}
\bpages{185--201}.
\bid{issn={0132-2818}, mr={0684472}}
\end{barticle}
%

\bptok{imsref}%
\endbibitem

\bibitem{Taq74}
\begin{barticle}[mr]
\bauthor{\bsnm{Taqqu},~\bfnm{Murad~S.}\binits{M.S.}}
(\byear{1975}).
\btitle{Weak convergence to fractional {B}rownian motion and to the {R}osenblatt process}.
\bjournal{Z. Wahrsch. Verw. Gebiete}
\bvolume{31}
\bpages{287--302}.
\bid{mr={0400329}}
\bptnote{check year}%
\end{barticle}
%

\bptok{imsref}%
\endbibitem

\bibitem{WanCai10}
\begin{barticle}[mr]
\bauthor{\bsnm{Wang},~\bfnm{Lihong}\binits{L.}} \AND
\bauthor{\bsnm{Cai},~\bfnm{Haiyan}\binits{H.}}
(\byear{2010}).
\btitle{Asymptotic properties of nonparametric regression for long memory random fields}.
\bjournal{J. Statist. Plann. Inference}
\bvolume{140}
\bpages{837--850}.
\bid{doi={10.1016/j.jspi.2009.09.012}, issn={0378-3758}, mr={2558409}}
\end{barticle}
%

\bptok{imsref}%
\endbibitem

\bibitem{Wan13}
\begin{barticle}[auto:parserefs-M02]
\bauthor{\bsnm{Wang},~\bfnm{Y.}\binits{Y.}}
(\byear{2014}).
\btitle{An invariance principle for fractional Brownian sheets}.
\bjournal{J. Theoret. Probab.}
\bnote{To appear}.
\bid{doi={10.1007/s10959-013-0483-2}}
\end{barticle}
%

\bptok{imsref}%
\endbibitem

\bibitem{Whi54}
\begin{barticle}[mr]
\bauthor{\bsnm{Whittle},~\bfnm{P.}\binits{P.}}
(\byear{1954}).
\btitle{On stationary processes in the plane}.
\bjournal{Biometrika}
\bvolume{41}
\bpages{434--449}.
\bid{issn={0006-3444}, mr={0067450}}
\end{barticle}
%

\bptok{imsref}%
\endbibitem

\end{thebibliography}
\end{document}